\documentclass[12pt]{amsart}
\usepackage{color}

\ifx\pdfoutput\undefined
\usepackage{graphicx}
\DeclareGraphicsExtensions{.pstex,.eps}

\else
\usepackage{amsfonts}
\usepackage[pdftex]{graphicx}  \DeclareGraphicsExtensions{.pdf,.mps}

\fi

\usepackage{amsxtra}
\usepackage{amssymb}
\usepackage{amsfonts}
\usepackage{amsmath}

\renewcommand{\Im}{\mathop{\rm Im}\nolimits}
\def\S{\mathhexbox278}

\newcommand{\beq}{\begin{equation}}
\newcommand{\ee}{\end{equation}}

\theoremstyle{plain} \newtheorem{theorem}{Theorem}[section]
\newtheorem{lemma}[theorem]{Lemma}

 \theoremstyle{definition}
 \theoremstyle{remark}
\newtheorem{remark}[theorem]{Remark}

\newcommand{\R}{{\mathbb R}}

\newcommand{\p}{{\partial}}

\def\im{{\rm i}}

\newcommand{\C}{\mathbb{C}}

\def\uno{{\kern+.3em {\rm 1} \kern -.22em {\rm l}}}

\numberwithin{equation}{section}

\begin{document}

\title
[Instability of approximate periodic solutions  for    NLS]
{On instability of some approximate periodic solutions  for the full
  nonlinear Schr\"odinger equation}

\author {Scipio Cuccagna and Jeremy L. Marzuola}

\begin{abstract}
Using the  Fermi Golden Rule analysis developed in \cite{cuccagnamizumachi}, we prove
asymptotic stability of asymmetric nonlinear bound states bifurcating
from linear bound states for a quintic nonlinear Schr\"odinger operator with
symmetric potential. This goes in the direction of proving that the approximate
periodic solutions  of the NLS in \cite{Marzuola} do not persist for the full NLS.
\end{abstract}

\maketitle

\section{Introduction}
\label{section:introduction}

We consider the quintic nonlinear
Schr\"odinger equation (NLS):
\begin{equation}
\label{NLS}
 \im u_{t }=- \partial _x^2 u +Vu -|u|^4 u , \, u(0,x)=u_0(x), \, (t,x)\in\mathbb{ R}\times
 \mathbb{ R}.
\end{equation}
We assume the following hypotheses.

\begin{itemize}
\item[(H1)]  The discrete spectrum of $-\partial_x^2 +V$ is $\sigma _d
  (- \partial _x^2  +V)= \{-\omega _0, -\omega _1\}$ such that $\omega
  _0-\omega _1$ is as small as we want and $0$ is not a resonance, i.e. if $u\in L^\infty (\R )$ satisfies
    $u''=Vu$, we have $u=0$.

 \item[(H2)] The potential $V(x)$ is even: $V(x)=V(-x)$ for all $x$ and $V(x)$ is real valued.

\item[(H3)] The potential $V$ is smooth and $\langle x \rangle ^n V^{(k)} (x) \in L^\infty (\R )$, for any $n$ and any $k$, where  $\langle x \rangle = \sqrt{1+x^2} .$

\item[(H4)]    Let $  \psi _j$ be real valued generators of $\ker (- \partial _x^2 +V +\omega _j)$  with $\| \psi _j \| _{L^2}=1$.   For $\langle f,  g\rangle =\int _\R f(x) g(x) dx$ we assume that for a fixed $a_0>0$
    \begin{equation}
  \label{eq:interaction} \begin{aligned} &5 \langle \psi _0^4,  \psi
    _1 ^2\rangle - \langle \psi _0^6,  1\rangle > a_0>0 \, , \end{aligned}
\end{equation}
  \begin{equation}
  \label{eq:interaction2} \begin{aligned} &
  5 \langle \psi _0^4,  \psi _1 ^2\rangle ^2 - \langle
     \psi _0^6,  1\rangle \langle \psi _0^2,  \psi _1 ^4 \rangle > a_0 \, .
 \end{aligned}
\end{equation}

\item[(H5)]  We assume that Fermi golden rule hypothesis (H5) stated in \S \ref{sec:FGR}.

\end{itemize}

\begin{remark}
\label{rem:regularity} The very strong regularity and decay hypotheses on
the potential $V(x)$  in (H3) are certainly  unnecessary.   See for example the dispersive estimates for for $e^{\im t (-\partial _x^2+V)}$
with
$V\in L^{1,1} (\R)$ in \cite{Goldberg}, or the case with
delta functions  in \cite{DMW}. Nonetheless, we do not  try to prove systematically
the estimates stated later in \S \ref{sec:dispersion} for these less regular potentials.
\end{remark}

We refer to the Appendix in \cite{kirr}  about the existence of
double well potentials  satisfying (H1)--(H5), as well as to the brief
computational discussion in Appendix \ref{sec:num}. Specifically, if one starts
with an even potential $V_0(x)$ such that  $- \partial _x^2  +V_0$
admits exactly one eigenvalue $-\Omega$, then setting 
$$V(x)=V_L(x)=V_0(x-L)+V_0(x+L)$$ for
$L\gg 1$ yields potentials with two eigenvalues, both very close to $-\Omega$.
Furthermore, for $\varphi (x)$ a normalized ground state for
$- \partial _x^2  +V_0$, then 
$$\psi _j (x)\approx 2^{-\frac{1}{2}}
(\varphi (x-L)+(-)^j\varphi(x+L) )  .$$ Then, 
$$\langle \psi _0 ^2, \psi
_1 ^4\rangle \approx \langle \psi _0 ^6, 1\rangle  \approx \langle
\psi _0 ^4, \psi _1 ^2\rangle  $$ 
because they are all about
$2^{\frac{11}{12}} \| \varphi \| ^{6}_{L^6}.$  Originally, results of
this nature appeared in the work of E. Harrell \cite{Har}.

  The  equations  \eqref{NLS} with $ V(x)=V_L(x)$  for $L\gg 1$   are the focus of recent research \cite{kirr,Marzuola}
  because of the rich patterns detected at small energies.
  The references \cite{kirr,Marzuola} both treat \eqref{NLS} with $|u|^4u$ replaced by a
  cubic nonlinearity (which for our purposes is a difficult
  problem due to the subcritical nature of the nonlinearity).  The
  result in \cite{kirr} proves the existence of a family of nonlinear
  ground states  $e^{\im t \omega }\phi _ {\omega }(x)$
of \eqref{NLS} that bifurcate out of the linear ground state
$e^{\im t \omega _0}\psi _ {0 }(x)$.  For $\omega > \omega _0$ close
to $ \omega _0$ the $e^{\im t \omega }\phi _ {\omega }(x)$ are
orbitally stable, see \S \ref{subsec:stability}, and even in $x$. At some
critical $\omega = \omega ^*$ the ground states bifurcate for $\omega > \omega ^*$  in two
families, one formed by even functions, which are unstable, and the
other formed by non symmetric functions, which are stable. The
arguments in \cite{kirr} are quite general, but here we double check
  that this behavior continues to hold also in the case of the
quintic NLS \eqref{NLS}.

In \cite{Marzuola}, the existence of
more complex long time patterns is analyzed  by studying the dynamics of a
simplified finite dimensional system, obtained by selecting a finite
number of variables of the NLS in an appropriate system of
coordinates. Over long times it is shown to be a good
approximation of the full NLS. In particular, this finite dimensional  approximation  of the NLS  admits a larger class of time periodic solutions
than just the standing waves. The question then becomes whether or
not these new periodic solutions persist also for the full NLS equation.
  In \cite{Marzuola} it is conjectured they do not persist. In this
paper we do not address the solutions considered in \cite{Marzuola},
but nonetheless for an easier problem we provide the mechanism by which the full NLS
disrupts periodic solutions of a simplified system similar to that in
\cite{Marzuola}.  In Appendix \ref{sec:fd} however, we will present
evidence that indeed similar dynamical solutions exist that would
collapse via the asymptotic stability analysis presented here to a
nonlinear bound state asymptotically.  Such dynamics have been abstractly studied in \cite{KirrKevPel} as well.

We recall that \cite{Marzuola} simplifies the NLS by
first choosing as system of coordinates the spectral decomposition of
$-\partial ^2_x+V$, and by then setting equal to $0$ the continuous components.
Here, we consider instead a natural representation of the portion of $H^1 (\R )$ near the surface of asymmetric ground states. There are then natural
finite dimensional approximations of the NLS admitting periodic solutions.
They are as legitimate approximate solutions of the NLS as those in
\cite{Marzuola}, although here we do not try to check as in \cite{Marzuola}
if they are good approximate solutions. Our solutions are relatively
easy because they live arbitrarily close to the surface of asymmetric
ground states.

When the full NLS \eqref{NLS} is turned on, these approximate periodic
solutions do not persist because the ground states are asymptotically stable.
Hence the periodic solutions of the simplified system, now split into a part
converging  in $H^1(\R )$  to the orbit of a ground state, and another
part which scatters like free radiation (see Theorem \ref{thm:as stab}).
This part of the paper fits easily in the framework of the literature
of asymptotic stability of ground states initiated in  \cite{SW1,SW2,BP1,BP2}.
We recall that the most general results are in \cite{Cu1}, which contains
a quite general proof of the so called Fermi golden rule. In the present
paper  though, we treat
a quite special situation, due to the hypothesis that $\sigma _d (- \partial _x^2  +V)$ consists of just two eigenvalues, and so it is enough to use
the simpler framework  of
\cite{cuccagnamizumachi,Cu5} (we recall that \cite{Cu5} is a revision and a simplification
of \cite{Cu4}, which contains various mistakes).   To address the solutions
in \cite{Marzuola} one can probably proceed similarly, though the
complexity of the dynamical systems studied and the cubic nonlinearity makes the analysis rather challenging. The difficulty
though, is that the solutions in \cite{Marzuola}, while of arbitrarily small
energy, might nonetheless not be sufficiently close to ground states.
 The issues then seem in some sense more "global", and closer in spirit
 to the problems addressed in
\cite{TY3,SW4}.

Following \cite{kirr}, we consider nonlinear ground states of the form
\begin{eqnarray*}
e^{i \omega t} ( c_0 \psi_0 + c_1 \psi_1 + \eta (x,t)) .
\end{eqnarray*}
Applying the machinery of \cite{kirr}, we prove that there is an $\omega ^*>
\omega _0 $, with $ \omega ^* -\omega _0\approx \omega _0  -\omega
_1$, such that for $\omega \in ( \omega _0, \omega ^*]$ there is a
uniquely defined map $\omega \to  \phi _ {\omega }$ in
 $C^{0}(( \omega _0,  \omega ^*], H^{k,s}(\R ))$ for any $(k,s)$,
 where
 \begin{equation}\label{eq:SobSp}
 \| u \| _{H^{k,s}}:= \| \langle x \rangle ^s (1-\partial _x^2)^{\frac{k}{2}} u \| _{L^{2}} \, .\end{equation}
  For $\omega > \omega ^*$ there is a bifurcation, with a branch of even ground states, and a branch of asymmetric ground states.
  We focus on the latter ones.  We then prove the following

\begin{theorem}\label{thm:as stab} There is a  $ \delta _0>0$ such that
 for any $ \omega \in (  \omega ^*,  \omega ^* + \delta _0 )$ there exist
  an $\epsilon_0>0$
and a $C>0$ such that if
$$ \|u_0-e^{\im \gamma}\phi_
{\omega  } \|_{H^1}<\epsilon <\epsilon _0,$$
there exist $\omega
_\pm \in ( \omega ^*,\omega ^* + \delta _0 )$, $\theta\in C^1(\R;\R)$ and $h _\pm \in H^1$
with
$$\| h_\pm\| _{H^1}+|\omega _\pm -\omega_1|\le C \epsilon $$
such that
\begin{equation}\label{scattering}
\lim_{t\to  \pm\infty}\|u(t,\cdot)-e^{\im \theta(t)}\phi_{\omega
_\pm}-e^{\im t\partial _x^2 }h _\pm\|_{H^1}=0 .
\end{equation}
It is possible to write
$$u(t,x)=e^{\im \theta(t)}\phi_{\omega  (t)}
+ A(t,x)+\widetilde{u}(t,x)$$
with $|A(t,x)|\le C_N(t) \langle x \rangle ^{-N}$ for any $N$, with $\lim _{|t|\to \infty }C_N(t)=0$,
with $\lim _{ t \to \pm \infty } \omega  (t)= \omega _\pm$, and such
that
the following Strichartz estimates are satisfied:
\begin{equation}\label{Strichartz} \|  \widetilde{u} \|
_{ L^\infty  _t( \mathbb{R},H ^{1 }_x(\mathbb{R})) \cap L^5 _t( \mathbb{R},W ^{1,10}_x(\mathbb{R}))
 \cap L^4_t( \mathbb{R},L^\infty _x(\mathbb{R}))}\le
 C\epsilon.
\end{equation}
 \end{theorem}

 Once the necessary spectral hypotheses in \cite{cuccagnamizumachi,Cu5} are proved in Section \ref{sec:discspec}, Theorem \ref{thm:as stab} is a direct consequence
 of \cite{cuccagnamizumachi,Cu5}. Nonetheless we give a sketch of the main
 steps in the proof. In particular we review in Section
 \ref{sec:dispersion}  the material on dispersion of linear operators
 needed in the proof. Here we recall that the absence of the endpoint
 Strichartz estimate on $\mathbb{R}$ requires
 some surrogates. The surrogates were found by Mizumachi \cite{M1}. However
 it turns out that \cite{M1} can be substantially simplified, and that
 the smoothing estimates contained in \cite{M1}, while interesting
 per se, are not necessary in the proof of the main result in \cite{M1}.
 In fact the classical smoothing estimates introduced by Kato in \cite{kato}
 are sufficient. This is discussed in \cite{cuccagnatarulli,Cu5} and
 is reviewed in Section \ref{sec:dispersion}.
  See also the recent results of
 \cite{DMW} to allow singular potentials in our analysis with
 restrictions to $k \leq 1$.

\section{Ground states for $\omega $ starting at $\omega _0$}
\label{sec:ground}

\subsection{Ground states for $\omega $ starting at $\omega _0$}
\label{subsec:near}

As in \cite{kirr} we consider ground states of the form
$$\phi _\omega =\rho _0\psi _0+\rho _1\psi _1+\eta (\psi _0,\psi _1, \omega),$$
with $\rho _0$ and $\rho _1$ in $\R$ and for $\eta (\psi _0,\psi _1, \omega)$
a real valued function belonging to $H^{k,s}(\R , \R )$ for any $(k,s)$
with $\langle \eta , \psi _j \rangle =0 $ for $j=0,1$.  We are looking
for the simplest asymmetric ground states possible and not for all
possible nonlinear ground states branching out of the linear ground state.

We denote by $P_c$ the
projection to the continuous spectral component of $-\partial _x^2+V$.
Hence we look at the system
\begin{equation}  \label{eq:LSreduction}\begin{aligned} & - \omega _0\rho _0+\omega \rho _0 -\langle \psi _0, (\rho _0\psi _0+\rho _1\psi _1+\eta )^5
\rangle =0\, , \\& - \omega _1\rho _1+\omega \rho _1 -\langle \psi _1, (\rho _0\psi _0+\rho _1\psi _1+\eta )^5\rangle  =0\, , \\& (-\partial _x^2+V+\omega )\eta = P_c(\rho _0\psi _0+\rho _1\psi _1+\eta )^5\,
 .
\end{aligned}
\end{equation}
By an elementary application of the implicit function theorem one
obtains the following

\begin{lemma} \label{lem:eta} For $\rho _0$ and $\rho _1$
sufficiently small, the third  equation \eqref{eq:LSreduction}
 admits a unique solution  $\eta = \eta (\rho _0,\rho _1, \omega)$ which depends smoothly in $(\rho _0,\rho _1, \omega)$
with values in $H^{k,s}(\R , \R )$ for any $(k,s)$ and
can be expressed as
\begin{equation}  \label{eq:eta} \eta (\rho _0,\rho _1, \omega ) =\sum
  _{j=0}^{4}\rho _0^{5-j}\rho _1^{ j} \eta _j (\rho _0 , \omega) +
  \rho _1^{ 5} E (\rho _0,\rho _1, \omega) ,
\end{equation}
where  $ \eta _j (\rho _0 , \omega)(x) =(-1)^j  \eta _j (\rho _0 , \omega)(-x) $ for all $x$.
\end{lemma}
The proof of Lemma \ref{lem:eta} is a standard application of the
impicit function theorem and a
resolvent identity.  Similar expansions are
proven in Propositions $4.1$ and $4.3$ in \cite{kirr}.

\begin{lemma} \label{lem:near} There is a fixed number $\varepsilon
  _0>0$ such that for $\rho _0\in [-\varepsilon _0, \varepsilon _0]$
  admits a unique   function $\omega =\omega (\rho _0)$, with $\omega
  (0)=\omega _0$, such that
$$( \rho _0, 0, \eta (\rho _0,0, \omega (\rho _0))$$
is a solution of system \eqref{eq:LSreduction}  in $C^{\infty} ((-\varepsilon _0,  \varepsilon _0), H^{k,s}(\R , \R ))$ for any $ (k,s)$.
 \end{lemma}
\proof  For  $\rho _1=0$ one can see that the third term in the lhs of
the second equation in \eqref{eq:LSreduction} is $0$, because it is
$\int _\R \psi _1  (\rho _0\psi _0 +\rho _0^5\eta _0)^5 dx$, which
vanishes since the integrand is an odd function. So the second equation  in \eqref{eq:LSreduction}  is trivial. Substituting  $\eta (\rho _0,0, \omega  )= \rho _0^5\eta _0(\rho _0 , \omega  )$ in the first, and factoring out a common factor $\rho _0$, we get
\begin{equation}
- \omega _0 +\omega   -\rho _0^4\langle \psi _0, (\psi _0+\rho _1\psi _1+\rho _0^4\eta _0 (\rho _0 , \omega  ))^5
\rangle =0.\nonumber
\end{equation}
By the implicit function theorem we can solve with respect to $\omega $
getting \begin{equation}\label{eq:omega near}
 \omega = \omega (\rho _0) = \omega_0   +\rho _0^4\langle \psi _0^6,
 1\rangle +O(\rho _0^8) .
\end{equation}
\qed

\begin{lemma} \label{lem:nearbranch} There is a number $\rho _0^*\in (0, \varepsilon _0)$, with $\rho _0^*\approx (\omega _0-\omega _1)^{\frac{1}{4}}$
such that at $\rho _0=\rho _0^*$ the function of Lemma \ref{lem:near}
satisfies also the equation
\begin{equation}  \label{eq:branch}\begin{aligned} &   - \omega _1 +\omega (\rho _0)  -\left \langle \psi _1, \frac{(\rho _0\psi _0+\rho _1\psi _1+\eta (\rho _0,\rho _1, \omega (\rho _0) )^5}{\rho _1}\big | _{ \rho _1 =0}\right \rangle  =0 \,
 .
\end{aligned}
\end{equation}
  \end{lemma}
\proof Equation \eqref{eq:branch} is, for $\eta _j=
\eta _j(\rho _0 , \omega (\rho _0))$, equivalent to
\begin{equation}  \label{eq:branch1}\begin{aligned} &   - \omega _1 +\omega (\rho _0)  -
\frac{\partial}{\partial \rho _1}\big | _{ \rho _1 =0}\left \langle \psi _1, (\rho _0\psi _0+\rho _1\psi _1+\eta (\rho _0,\rho _1, \omega (\rho _0) )^5 \right \rangle   \\& = - \omega _1 +\omega (\rho _0)  -5\rho _0^4
 \left \langle \psi _1, (\psi _0 +\rho _0^4\eta _0    )^4  (\psi _1 + \rho _0^4\eta _1)\right \rangle  =  \\& - \omega _1 +\omega _0 -    \rho _0^4 \left ( 5 \langle \psi _0^4,  \psi _1 ^2\rangle - \langle \psi _0^6,  1\rangle \right  ) + O(\rho _0^8)=0,
\end{aligned}
\end{equation}
where we have used \eqref{eq:omega near}.
By the implicit function theorem the last equation has exactly one solution
 \begin{equation}  \label{eq:branch2}  (\rho _0 ^*)^4= \frac{  \omega _0 -\omega _1}{5 \langle \psi _0^4,  \psi _1 ^2\rangle - \langle \psi _0^6,  1\rangle}+O((\omega _0 -\omega _1)^2).
 \end{equation}
\qed

At $\rho _0^*$ and at the corresponding value    $ \omega ^*=\omega (\rho _0^*)$, the family   of even ground states which we have found above, bifurcates in two families, one formed by even ground states and the other by asymmetric ground states. In the context of the cubic NLS, \cite{kirr} proves that for
$\omega >   \omega ^* $   the  even ground states are unstable while the
 asymmetric ground states are orbitally stable. In the rest of \S 2 we double check that asymmetric ground states have the same behavior of \cite{kirr}
 for our quintic NLS.

\subsection{Asymmetric ground states for $\omega >  \omega ^*$}
\label{subsec:asymm}
  Following \cite{kirr} we consider   the branch of
 asymmetric ground states defined for $\omega >  \omega ^*.$ Set
 \begin{equation}  \label{eq:asymm}\begin{aligned} & F(\rho _0,\rho _1, \omega  )= - \omega _0\rho _0+\omega \rho _0 -\langle \psi _0, (\rho _0\psi _0+\rho _1\psi _1+\eta (\rho _0,\rho _1, \omega   ) )^5
\rangle  \, , \\&  G(\rho _0,\rho _1, \omega  )= - \omega _1 +\omega  -\frac 1{\rho _1} \left \langle \psi _1,  (\rho _0\psi
  _0+\rho _1\psi _1+\eta (\rho _0,\rho _1, \omega  ) )^5  \right
\rangle   \,  .
\end{aligned}
\end{equation}
 We know that $F(\rho _0^*, 0,\omega ^* )=G(\rho _0^*, 0,\omega ^* )=0$.
 Hence, we apply the implicit function theorem and prove the following

  \begin{lemma} \label{lem:asymm} The Jacobian $\frac{\partial(F,G)}{\partial (\rho _0,\omega )}$ has rank 2 at $(\rho _0^*, 0,\omega ^* )$.   Correspondingly, by the implicit function theorem there are smooth functions
  $\rho _0(\rho _1)$ and $\omega (\rho _1)$ such that
  $$F(\rho _0(\rho _1), \rho _1,\omega (\rho _1) )=G(\rho _0(\rho _1),
  \rho _1,\omega (\rho _1) )=0,$$ which are defined for $\rho _1$ in a
  small neighborhood of 0 and such that
\begin{equation}  \label{eq:asymm1}\begin{aligned} & \rho _0(\rho _1) = \rho _0^* + \frac{\rho _0''(0)}{2}  \rho _1^2+o( \rho _1^2) , \\& \omega (\rho _1) = \omega^* + \frac{\omega''(0)}{2}  \rho _1^2+o( \rho _1^2) .
\end{aligned}
\end{equation}
We have \begin{equation}
 \label{eq:asymm2}\begin{aligned} & \omega''(0)= 20 (\rho
   _0^*)^2\frac{  5 \langle \psi _0^4,  \psi _1 ^2\rangle ^2 - \langle
     \psi _0^6,  1\rangle \langle \psi _0^2,  \psi _1 ^4 \rangle  }{5
     \langle \psi _0^4,  \psi _1 ^2\rangle - \langle \psi _0^6,
     1\rangle }+O((\rho _0^*)^8) > 0,\\
& \rho _0''(0)=\frac{5}{ \rho _0^* }  \frac{ \langle \psi _0^4, \psi _1^2 \rangle - \langle \psi _0^2, \psi _1^4\rangle  }{5 \langle \psi _0^4,  \psi _1 ^2\rangle - \langle \psi _0^6,  1\rangle } +O((\rho _0^*)^5) .
\end{aligned}
\end{equation}
\end{lemma}

 \proof We have
\begin{eqnarray*}
\partial _{\rho _0} F(\rho _0,\rho _1, \omega  )& = &   - \omega _0
+\omega   -5\langle \psi _0, (\rho _0\psi _0+\rho _1\psi _1 \\
&&+\eta (\rho _0,\rho _1, \omega   ) )^ 4 ( \psi _0+ \partial _{\rho _0}\eta (\rho _0,\rho _1, \omega   ) )
\rangle
 .
\end{eqnarray*}
For $\rho _1=0$ and $\omega =\omega (\rho _0)$ (see
 \eqref{eq:omega near}), we get
 \begin{equation}\label{Frho0}  \begin{aligned}  & \partial _{\rho _0} F=
  - \omega _0 +\omega   -5\rho _0^4\langle \psi _0, (\psi _0+\rho_0^4 \eta _0 )^ 4 ( \psi _0+ \rho _0^4\eta _1 )
\rangle \\
& = - \omega _0 +\omega   -5\rho _0^4 \langle \psi _0^6 ,1
\rangle +O(\rho _0^8)=   -4\rho _0^4 \langle \psi _0^6 ,1 \rangle
+O(\rho _0^8)  .
\end{aligned}
\end{equation}
 We have
\begin{equation}  \begin{aligned} \partial _{\omega} F=      \rho _0-5\langle \psi _0, (\rho _0\psi _0+\rho _1\psi _1+\eta (\rho _0,\rho _1, \omega   ) )^ 4   \partial _{\omega}\eta (\rho _0,\rho _1, \omega   ) )
\rangle  . \nonumber
\end{aligned}
\end{equation}
 For $\rho _1=0$, we get  \begin{equation} \label{Fomega} \begin{aligned} \partial _{\omega} F=      \rho _0-5 \rho _0^9 \langle \psi _0, (\psi _0+ \rho _0^4 \eta _0 )^ 4   \partial _{\omega}\eta _0 )\rangle =      \rho _0+O(\rho _0^9).
\end{aligned}
\end{equation}
 For $\rho _1=0$ we already know from  \eqref{eq:branch1} that
\begin{equation}   \begin{aligned} & G= - \omega _1 +\omega -5\rho _0^4
 \left \langle \psi _1, (\psi _0 +\rho _0^4\eta _0    )^4  (\psi _1 + \rho _0^4\eta _1)\right \rangle   .
\end{aligned} \nonumber
\end{equation}
 So,
\begin{equation}  \label{Grho0omega}\begin{aligned} & \partial _{\rho _0}G=   -20\rho _0^3
  \langle \psi _0^4,  \psi _1 ^2  \rangle  +O(\rho _0^7) \, ,  \quad  \partial _{\omega}G=   1  +O(\rho _0^8).
\end{aligned}
\end{equation}
 Then the Jacobian matrix  $\frac{\partial(F,G)}{\partial (\rho _0,\omega )}$ at $  (\rho _0, 0, \omega (\rho _0))$ is
 \begin{equation}  \label{Jacobian}\begin{aligned} &
 \begin{pmatrix}     -4\rho _0^4 \langle \psi _0^6 ,1 \rangle
+O(\rho _0^8) &
 \rho _0+O(\rho _0^9)   \\
-20\rho _0^3
  \langle \psi _0^4,  \psi _1 ^2  \rangle  +O(\rho _0^7)  &  1  +O(\rho _0^8)
 \end{pmatrix}\, ,
\end{aligned}
\end{equation}
 with inverse matrix
 \begin{equation}  \label{inverse}\begin{aligned} & M =   \frac{(4\rho
       _0^4)^{-1}}{  5 \langle \psi _0^4,  \psi _1 ^2\rangle - \langle
       \psi _0^6,  1\rangle +O(\rho _0^4)   } \times \\
&
 \begin{pmatrix}    1  +O(\rho _0^8)   &
 -\rho _0+O(\rho _0^9)   \\
 20\rho _0^3
  \langle \psi _0^4,  \psi _1 ^2  \rangle  +O(\rho _0^7)  &   -4\rho _0^4 \langle \psi _0^6 ,1 \rangle
+O(\rho _0^8)
 \end{pmatrix}.
\end{aligned}
\end{equation}
 We see below that $  \partial _{\rho _1} F= \partial _{\rho _1} G=0 $
 at $\rho _1=0$. This implies  immediately $\omega' (0)=\rho _0' (0)=0.$
 We have    \begin{equation}  \begin{aligned} &
 \partial _{\rho _1} F \big |_{\rho _1=0} =
 -5\langle \psi _0, (\rho _0\psi _0+\rho _1\psi _1+\eta   )^ 4
  ( \psi _1+ \partial _{\rho _1}\eta  )
\rangle \big |_{\rho _1=0}\\ & =  -5\langle \psi _0, (\rho _0\psi _0+\rho _0^5  \eta   )^ 4 ( \psi _1+ \rho _0^4 \eta _1 )
\rangle
=0 , \nonumber
\end{aligned}
\end{equation}
where the last equality is due to the fact that $\psi _1+ \rho _0^4 \eta _1$
is odd and the other factors are even. We have
 \begin{equation}  \begin{aligned} & \partial _{\rho _1}^2 F \big
     |_{\rho _1=0} = \\
&  -5\langle \psi _0, 4(\rho _0\psi _0+\rho _1\psi _1+\eta   )^ 3 (
\psi _1+ \partial _{\rho _1}\eta  )^2 \\
& +
 (\rho _0\psi _0+\rho _1\psi _1+\eta   )^ 4   \partial _{\rho _1}^2\eta
 \rangle \big |_{\rho _1=0}\\
& =
  -5\langle \psi _0, 4(\rho _0\psi _0 +\rho _0 ^5\eta  _0 )^ 3
  ( \psi _1+ \rho _0 ^4\eta  _1  )^2+
 (\rho _0\psi _0+\rho _0 ^5\eta  _0   )^ 4  2\rho _0 ^3  \eta _2
 \rangle  \\&
= -20 \rho _0 ^3\langle \psi _0^4, \psi _1^2 \rangle  +O(\rho _0 ^7) . \nonumber
\end{aligned}\nonumber
\end{equation}
 To show that $\partial _{\rho _1}G \big |_{\rho _1=0}=0$
 and compute  $\partial _{\rho _1}^2 G \big |_{\rho _1=0}$, we use
 the elementary fact that for two smooth functions $f(x)=xg(x)$ we have $g(0)=f'(0)$,
 $g'(0)=\frac{1}{2}f''(0)$ and $g''(0)=\frac{1}{3}f'''(0)$. Hence, calculating as above,
\begin{equation}  \begin{aligned} & \partial _{\rho _1}G \big |_{\rho _1=0}
  =  -\frac{1}{2}\partial _{\rho _1}^2 \big |_{\rho _1=0} \left \langle \psi _1,  (\rho _0\psi _0+\rho _1\psi _1+\eta  )^5  \right \rangle =\\&  -\frac{5}{2}  \left \langle \psi _1,  4(\rho _0\psi _0+\rho _0 ^5 \eta _0 )^3 ( \psi _1+ \rho _0 ^4 \eta _1  )^2 +(\rho _0\psi _0+\rho _0 ^5 \eta _0 )^4  2\rho _0 ^3 \eta _2  )
 \right \rangle =0,
 \end{aligned}\nonumber
\end{equation}
with the last equality due to the fact that we are integrating an odd
function.  Then,
\begin{equation}  \begin{aligned} & \partial _{\rho _1}^2G \big |_{\rho _1=0}
  =     \\
& -\frac{5}{3}\partial _{\rho _1}   \big |_{\rho _1=0}\left \langle
  \psi _1,  4(\rho _0\psi _0+\rho _1\psi _1 +  \eta    )^3 ( \psi
  _1+ \partial _{\rho _1} \eta    )^2 \right. \\
& \left. +(\rho _0\psi _0+\rho _1\psi _1 +  \eta )^4 \partial _{\rho _1}^2\eta
 \right \rangle  \\
& =  -\frac{5}{3}    \langle \psi _1,  12(\rho _0\psi _0+\rho _0 ^5
\eta _0 )^2  ( \psi _1+ \rho _0 ^4 \eta _1  )^3 +  (\rho _0\psi _0+\rho _0 ^5 \eta _0 )^4  6\rho _0 ^2 \eta _3 \\
& +   12(\rho _0\psi
_0+\rho _0 ^5 \eta _0 )^3 ( \psi _1+ \rho _0 ^4 \eta _1  )   2\rho _0
^3 \eta _2
   \rangle \\
& =-20 \rho _0 ^2 \langle \psi _0^2, \psi _1^4 \rangle +O(\rho _0 ^6) .
 \end{aligned}\nonumber
\end{equation}
Hence we get, using $M$ the matrix in \eqref{inverse} and the implicit differentiation formula
   $y' = - g_{y}^{-1}g_{x }$   for $g(x,y)=0$, and $y'' = -
   g_{y}^{-1}g_{x x }$ when $g_{x }=0$, $g_y$ invertible,
\begin{equation}  \begin{aligned} & \begin{pmatrix}    \rho _0''(0)
     \\
 \omega ''(0)
 \end{pmatrix} =  20\, M \begin{pmatrix}     \rho _0 ^3 \langle \psi _0^4, \psi _1^2 \rangle +O(\rho _0 ^7)
     \\
   \rho _0 ^2 \langle \psi _0^2, \psi _1^4 \rangle +O(\rho _0 ^6)
 \end{pmatrix} \\
& = \frac{ 5 }{  5 \langle \psi _0^4,  \psi _1 ^2\rangle
   - \langle \psi _0^6,  1\rangle +O(\rho _0^6)   } \\
& \times
 \begin{pmatrix}    \rho _0^{-1}  (\langle \psi _0^4, \psi _1^2 \rangle - \langle \psi _0^2, \psi _1^4\rangle )+ O(\rho _0 ^3)
     \\ 4\rho _0 ^2( 5 \langle \psi _0^4,  \psi _1 ^2\rangle ^2 - \langle \psi _0^6,  1\rangle  \langle \psi _0^2,  \psi _1^4\rangle  )+ O(\rho _0 ^{6})
 \end{pmatrix} .
 \end{aligned}\nonumber
\end{equation}

 \qed

\subsection{The stability of the asymmetric ground states} \label{subsec:stability}

We focus on  $$\phi    _{\omega  (\rho _1 )} =\rho _0(\rho _1 )\psi
_0+\rho _1\psi _1+\eta (\rho _0(\rho _1 ),\rho _1, \omega  (\rho _1 )
),$$
the asymmetric ground states.  Following \cite{kirr}, we prove now that
they are orbitally stable, that is for any such fixed $\omega =\omega  (\rho _1 )$ and
for any $\varepsilon >0$ there is $\delta >0$ such that for $u_0\in
H^1 (\R )$ such that
$$\sup _{\gamma \in \R}\|u_0-e^{\im \gamma}\phi_{\omega  } \|_{H^1}<\delta $$
and for $u(t)$ the solution of \eqref{NLS}
we have that $u(t)$ is globally defined and
$$\sup _{t,\gamma \in \R}\|u(t)-e^{\im \gamma}\phi_
{\omega  } \|_{H^1}<\varepsilon .$$
Set
$$ q(\omega ) =\|  \phi_{\omega  }  \| ^2 _{L^2} $$
and consider the pair of operators
\begin{equation}\label{eq:Lplus} \begin{aligned}&
L_{+}(\omega) = -\partial _x^2+V+\omega - 5 \phi _\omega ^4 \text{ and }
 L_{-}(\omega) = -\partial _x^2+V+\omega -   \phi _\omega ^4 .
\end{aligned}
\end{equation}
The proof of the orbital  stability of
asymmetric ground states is obtained in two steps. The second step is  the following well known result, see \cite{W1}, which in particular says
that the $\epsilon$ and the $\delta$ in the definition of orbital
stability can be here taken to be about of the same value. For the
proof, see \cite{Cu5}.

\begin{theorem}\label{thm:orbstab} Given the hypotheses and
conclusions of Lemma \ref{lem:stability} below, there $\exists \, \epsilon
_0>0$ and  $A_0(\omega )>0$ s.t. $ \epsilon \in (0,\epsilon _0)$
and
$$\| u(0,x) - \phi _ {\omega } \| _{H^1}<\epsilon $$
imply for the corresponding solution
$$ \inf  \{ \| u(t,x) -e^{i \gamma} \phi _ {\omega } (x ) \| _{H^1_x(  \R  )} : \gamma \in \R    \}  <
A_0(\omega ) \epsilon . $$
\end{theorem}

The first step to prove  the orbital  stability of
asymmetric ground states is the  following

\begin{lemma}
\label{lem:stability} Consider the asymmetric ground states
discussed in Subsection \ref{subsec:asymm} for $\omega >\omega ^*$. Then, for $\omega  $ sufficiently close to $ \omega ^*$
the following two statements are true.

 \begin{itemize}

\item [(1)] We have $\frac{d}{d\omega }q(\omega )>0$.

\item [(2)] For $\omega =\omega  (\rho _1 )$ the  set of eigenvalues of $L_{+}(\omega)$ is given by
$\sigma _d(L_{+}(\omega))=\{ -\mu _0  (\rho _1 ), \mu _1  (\rho _1 )
\}$, where
$$\mu _0  (\rho _1 ) \ge  4(\rho ^*_0)^4\langle \psi _0^6,
 1\rangle   + O((\rho ^*_0)^8)$$
and $\mu _1  (\rho _1 )>0$ with $\mu _1  (\rho _1 )\approx    (\rho
_0^*)^2 \rho _1^2   $.

\end{itemize}

\end{lemma}

\proof We have
\begin{equation}  \begin{aligned}& \frac{d}{d\rho _1 }q(\omega ) =
 \frac{d}{d\rho _1 } (\rho _0^2+\rho _1^2
 )+\frac{d}{d\rho _1 }\langle \eta , \eta \rangle
  = 2\rho _1  (1+ O(\rho _0)) \\& +\frac{d}{d\rho _1 } (  \rho
  _0^{10}\langle \eta _0 , \eta _0\rangle  +\rho _0^{8}\rho
  _1^{2}\langle \eta _1 , \eta _1\rangle  + 2 \rho _0^{8}\rho
  _1^{2}\langle \eta _0 , \eta _2 \rangle  +O(\rho _1^{3}) ) \\&
  =2\rho _1 (1+ O(\rho _0))  + \mathcal{O} (\rho_1^2).
\end{aligned}\nonumber
\end{equation}
We have $\frac{d}{d\rho _1 }\omega = \omega '' (0) \rho _1+O(\rho _1^2)$ where
$ \omega '' (0)>0$ by \eqref{eq:asymm2} and \eqref{eq:interaction2}. Claim $(1)$ follows from
\begin{equation}  \begin{aligned}&
\frac{d}{d\omega }q(\omega ) = \frac{d\rho _1}{d \omega} \,
\frac{d}{d\rho _1 }q(\omega )    = \frac{2\rho _1 (1 + O(\rho _0 ))}{\omega '' (0) \rho _1+O(\rho _1^2)} = \frac{2}{\omega '' (0)}(1 + O(\rho _0  ))>  0.
\end{aligned}\nonumber
\end{equation}
We consider the second claim.
First of all, by \eqref{eq:Lplus}  and by the fact that $\phi ^4_\omega$
is small, we know by general arguments that $\sigma _d(L_{+}(\omega))$ is
formed by two eigenvalues, the same number of
$\sigma _d( -\partial _x^2+V +\omega  )$ (recall also that in dimension 1 the eigenvalues have always multiplicity 1). Since
\begin{equation}  \begin{aligned} & \langle \psi _0,  L_+(\omega )  \psi _0
\rangle  = \langle \psi _0,  (-\partial _x^2+V  +\omega  )   \psi _0
\rangle  -5  \langle \psi _0^2,  \phi ^4_\omega \rangle  \\& =  \omega  - \omega _0 -5  \langle \psi _0^2,  (\rho _0 \psi _0+\rho _1\psi _1+\eta  ) ^4  \rangle
  .
\end{aligned}\nonumber
\end{equation}
At $  \omega =  \omega ^*$, $  \rho =  \rho _*$ and  $\rho _1=0$,
by \eqref{eq:omega near}
we obtain
\begin{eqnarray*}
\langle \psi _0,  L_+(\omega )  \psi _0
\rangle &  =&  \omega ^* - \omega _0 -5  \langle \psi _0^2,  (\rho
^*_0 \psi _0 +(\rho ^*_0)^5\eta _0 ) ^4 \\
& = & - 4(\rho ^*_0)^4\langle \psi _0^6,
 1\rangle   + O((\rho ^*_0)^8)<0  .
\end{eqnarray*}
This means that $L_+(\omega ^*)$ has at least one negative eigenvalue. By
continuity    for $\omega >\omega ^*$ close to $\omega ^*$, $L_+(\omega )$
has at least one negative eigenvalue, which we denote by  $-\mu _0 (\rho _1)$.
We now start the discussion on $ \mu _1 (\rho _1)$.
We claim  that
\begin{equation} \label{eq:L+crit} \begin{aligned} &L_+( \omega ^*)(\psi _1+(\rho _0^*)^4\eta _1(\rho _0^*,\omega ^* ) )  =0.
\end{aligned}
\end{equation} Starting by $L_+(\omega  )\partial _\omega \phi  _\omega =- \phi  _\omega $ and by Lemma \ref{lem:exp phipr} below, we get
\begin{equation}  \begin{aligned} &
L_+(\omega  )\partial _\omega \phi  _\omega =  \frac{L_+(\omega  )}{\omega ''(0)\rho _1 +O((\rho _1)^2)} (\psi _1+ \rho _0 ^4\eta _1(\rho _0 ,\omega   ) + O(\rho _1))=- \phi  _\omega.
\end{aligned}\nonumber
\end{equation}
Multiplying the equations by $\rho _1$ and letting $\rho _1\searrow 0$
we obtain \eqref{eq:L+crit}. This implies that for small $\rho _1$ the eigenvalue  problem $L_+(\omega )g=\mu g$ admits a small solution  $ \mu =\mu _1 (\rho _1)$
with $ \mu _1 (0) =0$. There is a unique
  solution $L_+(\omega )g=\mu g$ of the form $g=\alpha  _0 \psi _0+  \psi _1+\xi $, where   $\xi
=P_cg$,  $\mu  =\mu _1$ and $g$ are functions of $\rho _1$.
 The expression $L_+(\omega )g=\mu g$ is equivalently expressed as follows considering a Lyapunov Schmidt reduction as in \S 5 in \cite{kirr}:
\begin{equation}  \label{eq:L+}\begin{aligned} &
\langle \psi _0, L_+(\omega )(\alpha _0 \psi _0+  \psi _1+\xi)
\rangle  =  \mu \alpha _0 \, , \\& \langle \psi _1, L_+(\omega
)(\alpha _0 \psi _0+  \psi _1+\xi) \rangle  =  \mu  \, ,
  \\& (-\partial _x^2+V+\omega -\mu  )\xi = 5P_c \phi _\omega ^4   (\alpha _0 \psi _0+  \psi _1+\xi) \,
 .
\end{aligned}
\end{equation}
  We set
$\xi _*=(\rho _0^*)^4  \eta _1^*  ,$
 where $\eta _j^*=\eta _j(\rho _0^*, \omega ^* ) $. We have
\begin{equation}  \label{eq:alpha0}\begin{aligned} &
\left (  \omega - \omega _0-\mu -5\langle \phi _\omega ^4, \psi _0^2
\rangle \right ) \alpha _0   =5\langle \phi _\omega ^4, \psi _0 (\psi
_1+ \xi ) \rangle
 .
\end{aligned}
\end{equation}
By \eqref{eq:L+crit}  we know $\alpha _0  (\rho _1)_{| \rho
_1=0}=0.$ At $\rho _1=0 $ we have for $\omega = \omega ^*$
\begin{equation}  \label{eq:alpha01}\begin{aligned} &
\left (  \omega   -5\langle \phi  _{\omega  } ^4, \psi _0^2 \rangle
\right ) \alpha _0'   =20\langle \phi _\omega ^3 (\psi _1+\xi _* ),
\psi _0 (\psi _1 +\xi _*) \rangle +5\langle \phi _\omega ^4, \psi _0
\partial _{\rho _1}\xi \rangle
 .
\end{aligned}
\end{equation}
At $\rho _1=0 $, from \eqref{eq:L+} and  from $\omega '(0)=0$, see \eqref{eq:asymm1}, we get
\begin{equation}  \label{eq:derxi}  \begin{aligned} & \partial _{\rho _1}\xi =
5R_{V}(-\omega ^*)  \left [  4 \phi _{\omega ^*}^3 (\psi _1 + \xi
_*)^2+ \phi _{\omega ^*}^4 \alpha _0'\psi _0+ \phi _{\omega ^*}^4
\partial _{\rho _1}\xi  \right ]\\&  +\mu ' (0) R_V(-\omega ^*)\xi _*
\end{aligned}
\end{equation}
for $R_V(z)= (-\partial _x^2+V-z)^{-1}$.
We have
\begin{equation}  \label{eq:mu 0}\begin{aligned} &  \mu = \omega -
\omega _1 -5  \alpha _0 \langle \phi _\omega ^4, \psi _0 \psi _1
\rangle -5    \langle \phi _\omega ^4,   \psi _1 ( \psi _1+ \xi )\rangle
 .
\end{aligned}
\end{equation}
Then, \eqref{eq:derxi} and \eqref{eq:mu 0} imply  $\mu '(0)=0$.
Indeed, at $\rho _1=0$ we have $\omega '(0)=\alpha _0 (0)=0$, hence we see
\begin{equation}  \label{eq:mu 1}\begin{aligned} &  \mu '(0) =   -5  \alpha _0' \langle \phi _{\omega ^* } ^4, \psi _0 \psi _1
\rangle - 20    \langle \phi _{\omega ^* }^ 3,   \psi _1 ( \psi _1+ \xi _*)^2\rangle -\\&  5    \langle \phi _{\omega ^* } ^4,   \psi _1  \partial _{\rho _1} \xi  \rangle  = - 5    \langle \phi _{\omega ^* } ^4,   \psi _1  \partial _{\rho _1} \xi  \rangle  = - 5   \mu '(0)  \langle \phi _{\omega ^* } ^4,   \psi _1  R_V(-\omega ^*)\xi _*  \rangle .
\end{aligned}  \nonumber
\end{equation}
 Since in \eqref{eq:derxi} we have $\mu '(0)=0$ and by the fact that the resolvent
 $R_V(z)$ preserves the spaces of  even (resp. odd)  functions,
 we conclude that
  $\partial _{\rho _1}\xi _{| \rho _1=0}$ is even. We also have the estimate
  \begin{equation*} \| \partial _{\rho _1}\xi _{| \rho _1=0} \| _{L^\infty} \le C
(\rho _0^*) ^3 +C (\rho _0^*) ^4\alpha _0'(0).\end{equation*}
This and \eqref{eq:alpha01} yield
\begin{equation}  \label{eq:alpha02}\begin{aligned} &
\left (     -4 \langle   \psi _0^6,1 \rangle (\rho _0^*) ^4 +O((\rho
_0^*) ^7) \right ) \alpha _0'(0) =20  (\rho _0^*) ^3  \langle \psi
_0 ^4,\psi _1^2 \rangle +O((\rho _0^*) ^5)
\end{aligned}\nonumber
\end{equation}
and so   $\alpha _0'(0) =-4(\rho _0^*) ^{-1} \frac{\langle \psi _0
^4, \psi _1 ^2\rangle}{\langle \psi _0 ^6,1  \rangle} +O(\rho _0^*)$.
We have
\begin{equation}  \label{eq:mu 01}\begin{aligned} &  \mu ''(0)= \omega ''(0)
-40  \alpha _0'(0) \langle \phi _\omega ^3 (\psi _1 +\xi _*), \psi
_0 \psi _1 \rangle   \\& -20 \langle \phi _\omega ^3\partial _{\rho
_1}^2 \phi _\omega, \psi _1  (\psi _1+\xi _*  ) \rangle  -60 \langle \phi _\omega ^2,
(\psi _1+\xi _*)^3 \psi _1    \rangle   \\ &
-40 \langle \phi _\omega ^3(\psi _1+\xi _*)  , \psi _1
\partial _{\rho _1}\xi   \rangle
-5 \langle \phi _\omega ^4, \psi _1 \partial _{\rho _1}^2\xi\rangle
 .
\end{aligned}
\end{equation}
Proceeding as above,  $|\alpha _0''(0) | \le
 C(\rho _0^*) ^{-2} $ and $ \| \partial _{\rho _1}^2\xi _{| \rho
_1=0} \| _{L^\infty} \le C  .$ The terms in the second and
third lines of \eqref{eq:mu 01}
 are $O((\rho _0^*) ^{3})$. As a result,
\begin{equation}  \label{eq:mu 02}\begin{aligned}
     \mu ''(0)&=  20 (\rho _0^*)^2\frac{  5 \langle \psi _0^4,
  \psi _1 ^2\rangle ^2 - \langle \psi _0^6,  1\rangle \langle \psi _0^2,
    \psi _1 ^4 \rangle  }{5 \langle \psi _0^4,  \psi _1 ^2\rangle
    - \langle \psi _0^6,  1\rangle } +160(\rho _0^*)^2\frac{\langle \psi _0
^4, \psi _1 ^2\rangle ^2}{\langle \psi _0 ^6,1  \rangle}\\& -60(\rho
_0^*)^2\langle \psi _0 ^2, \psi _1 ^4\rangle   +O((\rho _0^*) ^{3})>0
\end{aligned}
\end{equation}
by \eqref{eq:interaction}--\eqref{eq:interaction2}.
\qed

 \begin{lemma}
  \label{lem:exp phipr}  At $\rho _1=0$ for the asymmetric branch we have the expansion
  \begin{equation}  \label{eq:phipr}\begin{aligned}&
  \partial _\omega \phi _ \omega = O(\rho _1)+ \frac{1}{\omega '(\rho _1 ) }  (\psi _1+(\rho _0^*)^4\eta _1(\rho _0^*,\omega ^* ) ) + \\&\frac{\rho _0 ''(0 )}{\omega ''(0 ) } (\psi _0+ 5(\rho _0^*)^4\eta _0(\rho _0^*,\omega ^* )) +
  (\rho _0^*)^5\partial _\omega \eta _0(\rho _0^*,\omega ^* ) .
  \end{aligned}
\end{equation}
\end{lemma}
\proof By Lemma \ref{lem:eta} we have the following formula
 which follows from\begin{equation}  \begin{aligned} &
 \partial _\omega \phi  _\omega = \frac{1}{\frac{d\omega }{d\rho _1}} \partial _{\rho _1}  (\rho _0(\rho _1 )\psi _0+\rho _1\psi _1+\eta (\rho _0(\rho _1 ),\rho _1, \omega  (\rho _1 ) ))\\& =  \frac 1{\omega ' } (\psi _1+ \rho _0 ^4\eta _1(\rho _0 ,\omega   )+ \rho _0 '(\psi _0 + 5 \rho _0^4 \eta _0)  +
  \omega ' \rho _0^5 \partial _\omega \eta _0+ O(\rho _1^2)) .
\end{aligned}\nonumber
\end{equation}
We the
use Lemmas \ref{lem:asymm}.
\qed

\section{The discrete spectrum of the linearization}
\label{sec:discspec}

Consider the operators
\begin{equation}  \label{eq:lineariz}\begin{aligned} & \mathcal{L}_{\omega}=
\begin{pmatrix}     0 & L_-(\omega )   \\ -L_+(\omega ) & 0
 \end{pmatrix} \text{ and } J=
\begin{pmatrix}     0 & 1   \\ -1 & 0
 \end{pmatrix} .
\end{aligned}
\end{equation}

 \begin{lemma}
  \label{lem:discspec}  We have $\sigma _d(\mathcal{L}_{ \omega ^*})={0}$. For $\{ \}$ meaning span, we have
  \begin{equation}  \label{eq:ker}
  \begin{aligned} &  \ker \mathcal{L}_{ \omega ^*} =\{  e_1(\omega ^*), e_2(\omega ^*) \} \text{ with } e_1(\omega  ) =
\begin{pmatrix}      0   \\  \phi _{\omega}
 \end{pmatrix}  \, , \,  e_2(\omega ^* ) =
\begin{pmatrix}      \beta    \\   0
 \end{pmatrix} ,
\end{aligned}
\end{equation}
where we set $\beta =\psi _1+(\rho _0^*)^4\eta _1^* $.
The generalized kernel is
\begin{equation}  \label{eq:genker}
  \begin{aligned} &  N_g( \mathcal{L}_{\omega ^*}) =\{  e_j(\omega ^*): j=1,...4 \}
\end{aligned}
\end{equation}
with
\begin{equation*}
e_3(\omega ^* ) =
\begin{pmatrix}     \alpha   \\   0
 \end{pmatrix}  \, , \,  e_4(\omega ^*  ) =
\begin{pmatrix}      0      \\   \gamma
 \end{pmatrix} ,
\end{equation*}
where  $  \phi _{\omega ^*} = L_+(\omega ^* )\alpha  $  and    $\beta   = L_-(\omega ^* )\gamma $ .  We have that $\alpha (x)$ is an even function, while $\beta (x)$ and $\gamma  (x)$ are odd functions.
\end{lemma}
 \proof  The relation \eqref{eq:ker} follows by the fact that $0$ is an eigenvalue of $L_-(\omega )$ and   $L_+(\omega ^* )$.
  The equations $  \phi _{\omega ^*} = L_+(\omega ^* )\alpha  $
  and  $\beta   = L_-(\omega ^* )\gamma  $
  admit    solutions because $\langle \phi _{\omega ^*} ,\beta \rangle =0 $, since
  $\phi _{\omega ^*}$ is even and $\beta$ odd. In fact $\alpha $ coincides with the second line of formula \eqref{eq:phipr}. We conclude that
  $N_g( \mathcal{L}_{\omega ^*})$ contains  the rhs of
 \eqref{eq:genker}   and so $\dim N_g( \mathcal{L}_{\omega ^*})\ge  4$.
 To see that they are equal it is enough to check that $\dim N_g( \mathcal{L}_{\omega ^*})\le  4$.   $\mathcal{L}_{\omega ^*}$ is a small perturbation of $J (-\partial _x^2+V +\omega ^* )$ which has $4$ eigenvalues, all close to 0. This implies $\dim N_g( \mathcal{L}_{\omega ^*})\le 4$.
 \qed

 \begin{lemma}
  \label{lem:discspec1}  Let $ \omega > \omega ^*$. Then
    $\sigma _d(\mathcal{L}_{\omega  })=\{0, \im \lambda (\omega) ,-\im \lambda (\omega) \}$
    with $\lambda (\omega) >0 $ a simple eigenvalue.
    \end{lemma}
     \proof  We know that $ \dim  N_g(\mathcal{L}_{\omega  })=2$, with $  N_g( \mathcal{L}_{\omega  })=\{  e_1(\omega  ),e_2(\omega  )  \}$ with $e_1(\omega  ) $ as in Lemma \ref{lem:discspec} and with
 $(e_2(\omega  ))^T =(\rho _1 \partial _{\omega}\phi _{\omega  } , 0   )$.  Let $D $ be a small disk containing
  the origin in its interior. We know that $\mathcal{L}_{\omega} -\mathcal{L}_{\omega ^*}$ is continuous in $\omega $ with values in the space of bounded operators from  $L^2(\R ) $ in itself with the uniform topology, and that  $\mathcal{L}_{\omega ^*}$ has exactly $4$ eigenvalues inside $D$
  (in fact just $1$ with algebraic multiplicity $4$) and that
  $\partial D $ is in the resolvent set of $\mathcal{L}_{\omega ^*}$. Then,
 $$\sigma _d(\mathcal{L}_{\omega  })\cap D=\{0, \im \lambda (\omega)
 ,-\im \lambda (\omega) \}$$
follows by the fact that $\sigma _d(\mathcal{L}_{\omega  })$ is symmetric with respect to the coordinate axes
 and $0$ has algebraic multiplicity $2$ for  $\mathcal{L}_{\omega}$, as we reminded above, and that the two eigenvalues cannot lie in $\R$ (since we know that
 the $\phi _\omega $ for $\omega >\omega ^*$ are orbitally stable).

 By standard arguments in perturbation theory, exploiting only
$$|\phi _{\mu} (x)|\le C e^{-a|x|}$$
for $a>0$ and $C>0$ fixed and for both $\mu =\omega ^*, \omega$, it is possible to prove
that
$$\sigma _d(\mathcal{L}_{\omega  })\cap (\C \backslash D)$$
is empty for $\omega $ close enough to $\omega ^* $. Specifically, and
quite informally,  elements of
$$\sigma _d(\mathcal{L}_{\omega  })\cap
(\C \backslash D)$$
 could originate by singularities of
$(\mathcal{L}_{\omega ^* }-z) ^{-1}$ on a second sheet of the Riemann
surface where it is defined, or by the points $\pm \im \omega ^*$ if 0
was a resonance of $-\partial _x^2+V$. But the latter is excluded  by
hypothesis and the former can be ruled out for $\omega$ close enough
to $\omega ^*$. We skip the details: the analysis at the endpoints
 is similar to material in \cite{Cu6}; the analysis of the eigenvalues coming from the second sheet can be derived from \cite{CPV} .
\qed

\begin{lemma}
  \label{lem:lambda}  Consider for $ \omega > \omega ^*$ the
  eigenvalue from Lemma \ref{lem:discspec1} with
  $\lambda (\omega) >0 $.
Let $ \omega > \omega ^*$, then there exists a fixed $C >0$ such
that
  $\lambda (\omega (\rho _1)) >C \rho _1\rho _0^*$.
    \end{lemma}
\proof We recall that if $\mathcal{L}_{\omega} U = \im \lambda U $
and   $U^T =(u,v)$, then
$$L_-(\omega   )L_+(\omega   ) u=\lambda ^2u.$$
Then, $\langle u , \phi _{\omega}\rangle =0$ and one can define $f =(L_-(\omega   )) ^{-\frac{1}{2}} u$ which is s.t.
$$(L_-(\omega   )) ^{ \frac{1}{2}}L_+(\omega   )(L_-(\omega   )) ^{
  \frac{1}{2}} f=\lambda ^2f.$$
Since  one can proceed backwards, we have
\begin{equation} \begin{aligned}
     &   \lambda ^2 =\min _{\langle g , \phi _{\omega}\rangle =0}\frac{ \langle (L_-(\omega   )) ^{ \frac{1}{2}}L_+(\omega   )(L_-(\omega   )) ^{ \frac{1}{2}}g ,  g\rangle}{\|  g \| ^2 _{L^2}} \ge \\&  \min _{\langle f , \phi _{\omega}\rangle =0}\frac{ \langle  L_+(\omega   )f, f\rangle}{\| f \| ^2 _{L^2}} \min _{\langle g , \phi _{\omega}\rangle =0}\frac{ \langle L_-(\omega   ) g ,  g\rangle}{\|  g \| ^2 _{L^2}}.
\end{aligned}\nonumber\end{equation}
We prove now that
\begin{equation}\label{itermid1} \begin{aligned}
     &     \min _{\langle f , \phi _{\omega}\rangle =0}\frac{ \langle  L_+(\omega   )f, f\rangle}{\| f \| ^2 _{L^2}} > C_1 (\rho _0^*)^2\rho _1^2 .
\end{aligned} \end{equation}
Let $\| f \| ^2 _{L^2}=1$ and $\phi = \phi _{\omega}$.
Let $ L_+(\omega   )\chi _j =(-) ^{j+1}\mu _j \chi _j$  for $j=0,1$ with
$\| \chi _j   \| _{L^2} =1$.
For a function $g$, let $g_j = \langle g , \chi _j\rangle $ and let
$g_c=\| P_c(L_+(\omega ))g \| _{L^2} $, where (only for this proof)
$$P_c g:=g- g_0  \chi _0 -g_1  \chi _1 .$$
Then,
$$\langle  L_+(\omega   )f, f\rangle \ge -\mu _0 f_0^2+\mu _1 f_1^2 +
\omega  f_c ^2 =:F(f_0 ,f_1,f_c).$$

We will prove then that, subject to the constraints in the last two lines of \eqref{itermid2} below,
we have  $F> C_2 \rho _1\rho _0^* .$ Since  $F$ is continuous
 for the strong and weak topology in $L^2$, there
 exists a constrained minimizer.
This implies that, for $a$ and $b$ Lagrange multipliers, we have:
\begin{equation}\label{itermid2} \begin{aligned}
     &      2(\omega -a)P_cf=b P_c\phi , \\  -&2\mu _0 f_0= 2a f_0 +b
     \phi _0 , \\ &  2\mu _1 f_1= 2a f_1+b \phi _1 , \\ & f_0^2+
     f_1^2 +    f_c ^2 =1 , \\ & f_0\phi_0+  f_1\phi_1 +   \langle P_c f,P_c \phi \rangle =0.
\end{aligned} \end{equation}
For  $P_c\phi  $  and  $P_cf  $ proportional to each other, the last equation in
\eqref{itermid2}  is the same as $f_0\phi_0+  f_1\phi_1 +   f_c\phi_c =0$.
If $P_c\phi$ and $P_cf$ are not proportional, then $b=0$ and $\omega =a$. Then $ (\omega +\mu _0) f_0=0$ implies $f_0=0$ since $\omega +\mu _0>0$.
Given $ (\omega -\mu _1) f_0=0$, we have $f_1=0$ since $\mu _1=O(\rho _1^2)$
while $\omega > \omega ^* >0$. Then, $f_c=1$ with $F=\omega$, which is clearly the maximum value, and not the minimum. Hence we can assume that $P_c\phi  $  and  $P_cf  $ are proportional. Then we minimize $F$ under the constraint
\begin{equation}\label{itermid} \begin{aligned}
     &    f_0^2+  f_1^2 +    f_c ^2 =1 , \\ & f_0\phi_0+  f_1\phi_1 +   f_c\phi_c =0.
\end{aligned} \end{equation}
Notice that the plane
 can be parametrized by $f_0= \phi _1 u + \phi _c v$,
   $f_1= -\phi _0 u  $,  $f_c= - \phi _0 v$. Then,
  \begin{equation}\label{itermid3} \begin{aligned}
     &  F(\phi _1 u + \phi _c v ,-\phi _0 u,- \phi _0 v) = -\mu _0 (\phi _1 u + \phi _c v) ^2+ \mu _1 \phi _0^2 u^2+ \omega \phi _0^2 v^2\\& = (-\mu _0 \phi _1^2+ \mu _1  \phi _0^2) u^2  + (-\mu _0 \phi _c^2+ \omega  \phi _0^2) v^2- 2 \mu _0 \phi _1\phi _c  uv .
\end{aligned} \end{equation}
This is a quadratic form in $(u,v)$ with eigenvalues, $x$, the roots of
 \begin{equation}\label{itermid31} \begin{aligned}
     &  (x- (-\mu _0 \phi _1^2+ \mu _1  \phi _0^2)) ( x-  (-\mu _0 \phi _c^2+ \omega  \phi _0^2)) - \mu _0^2 \phi _1^2\phi _c^2     \\&  = x^2 - (-\mu _0 \phi _1^2+ \mu _1  \phi _0^2-\mu _0 \phi _c^2+ \omega  \phi _0^2) \, x+\\& +(-\mu _0 \phi _1^2+ \mu _1  \phi _0^2) (-\mu _0 \phi _c^2+ \omega  \phi _0^2)   -  \mu _0^2 \phi _1^2\phi _c^2  =0.
\end{aligned} \end{equation}
We have   \begin{equation}\label{itermid4} \begin{aligned}& -\mu _0 \phi _1^2+ \mu _1  \phi _0^2-\mu _0 \phi _c^2+ \omega  \phi _0^2 = \omega _0 ( \rho _0^*)^2 +o(( \rho _0^*)^2)
\end{aligned} \end{equation}
and
\begin{equation}\label{itermid5} \begin{aligned}&          -\mu _0  \omega\phi _1^2   \phi _0^2-\mu _0 \mu _1 \phi _c^2 \phi _0^2 -  \mu _0^2 \phi _1^2\phi _c^2+ \mu _0^2  \phi _1^2 \phi _c^2+\mu _1\omega    \phi _0^4 \\& =  (\mu _1\omega \phi _0^2-\mu _0  \omega   \phi _1^2-\mu _0 \mu _1 \phi _c^2       )   \phi _0^2.
\end{aligned} \end{equation}
We claim that \eqref{itermid5}$\approx \omega _0(\rho _0^*)^6 \rho _1^2 . $
This and \eqref{itermid4} imply that the polynomial in \eqref{itermid31}
has both roots positive, one about $\omega _0(\rho _0^*)^2 $  and the other about $(\rho _0^*)^4 \rho _1^2 $. Then, the minimum of
\eqref{itermid3}  for $u^2+v^2=1$ is about $(\rho _0^*)^4 \rho _1^2 $.
Since $f_0^2+  f_1^2 +    f_c ^2 =1$  implies  $u^2+v^2\approx (\rho
_0^*)^{-2}$, it follows that the minimum of $F$ is $>C (\rho
_0^*)^2\rho _1^2 $ for a fixed $C$.

To show \eqref{itermid5}$\approx \omega _0(\rho _0^*)^4 \rho _1^2  $,
we observe that since $L_+(\omega ) =L_-(\omega )-4\phi _\omega ^4,$ we have
$\mu _0\le 4 \| \phi _\omega \| _\infty^4 \le C \rho _0^4$ for a fixed $C$.
Then, the claim follows from
 \begin{equation}  \begin{aligned}&  |\text{\eqref{itermid5}} |  \ge
     \phi _0^2(\omega\mu _1\phi _0^2-C \omega \rho _0^4  \rho _1^2- C
     \rho _0^4 \mu _1 \phi _c^2  )  \\
&>  \phi _0^2\mu _1(\frac{1}{2}\omega\phi _0^2 - C \rho _0^4   \phi _c^2  )   >  \frac{1}{3}\phi _0^4\mu _1\omega,
\end{aligned} \nonumber \end{equation}
where we exploit  $\phi _1=O(\rho _1)$, $\phi _0\approx \rho _0 $, $\phi _c=O(\rho _c)$ and $\mu _1\approx \rho _0^2  \rho _1^2\gg \rho _0^4  \rho _1^2. $ Then $|\text{\eqref{itermid5}} |\gtrsim  \rho _0^6 \rho _1^2. $ \qed

\section{Set up for Theorem \ref{thm:as stab} and dispersion for the linearization}
\label{sec:dispersion} Theorem \ref{thm:as stab} is a consequence of
\cite{Cu5}. Notice that since the linearization has just one pair of
nonzero eigenvalues, the hamiltonian set up in \cite{Cu1} is
unnecessary, and the theory in \cite{Cu5,cuccagnamizumachi} is
adequate. We recall that due to absence of the endpoint Strichartz
estimate in 1 D, the theory requires some adequate surrogate. This
for 1 D was provided by Mizumachi \cite{M1}. The theory in \cite{M1}
though, is more complicated then necessary. The simplifications
were provided in \cite{Cu5,cuccagnatarulli}. Subsequent papers  like \cite{KPS,PelinovskyStefanov} return to more complicated   approach \cite{M1}. So, even though Theorem \ref{thm:as stab} is a direct
consequence of \cite{Cu5}, we will take the opportunity to state the
various steps of the proof, in order also to point out the points in
\cite{M1} and in \cite{PelinovskyStefanov} which can be simplified.

Recall the ansatz \begin{equation}\label{eq:ansatz}  u(t,x) = e^{\im \Theta (t)} (\phi _{\omega (t)} (x)+
r(t,x)) \, , \, \Theta (t)= \int _0^t\omega (s) ds +\gamma (t). \end{equation} Inserting the ansatz into the NLS \eqref{NLS} we get
\begin{equation} \begin{aligned}
&
  \im  r_t  =
 -  r _{xx}+V r+\omega (t) r-
3\phi _{\omega (t)}
^4 r  -\phi _{\omega (t)}
^4  \overline{  r }\\& + \dot \gamma (t) \phi
_{\omega (t)} - \im \dot \omega (t)
\partial _\omega \phi   _{\omega (t)}
+ \dot \gamma (t) r
 +
  O(r^2).\end{aligned}\nonumber
\end{equation}
We set
$^tR= (r,\bar r) $,  $^t\Phi = ( \phi _{\omega } , \phi _{\omega } )
$ (using a different frame from the one in \S \ref{sec:discspec}
 and we rewrite the above equation as
\begin{equation}\label{eqR} \begin{aligned}
& \im   R _t =H _{\omega}   R +\sigma _3 \dot \gamma   R
+\sigma _3 \dot \gamma \Phi - \im  \dot \omega \partial _\omega \Phi
+O(R^2) \end{aligned}
\end{equation}
where:
 \begin{equation} \label{eq:operator}\begin{aligned} &\sigma _1=
\begin{pmatrix}0 &
1  \\
1 & 0
 \end{pmatrix} \, ,
\sigma _2=\begin{pmatrix}  0 &
\im  \\
-\im & 0
 \end{pmatrix} \, ,
\sigma _3=\begin{pmatrix} 1 & 0\\0 & -1 \end{pmatrix}
\, ; \\ &
H_{\omega,0}=\sigma_3(-\partial ^2_x+V+\omega),\,
\mathcal{V}_\omega=-3\sigma _3 \left[ \phi ^4_{\omega }  \right]
+2\im   \phi ^4_{\omega }\sigma _2;
\\ &   H _{\omega} =H_{\omega,0}+V_\omega .
\end{aligned}
\end{equation}
We know that $0$ is an isolated eigenvalue of $H_\omega$,
$\dim N_g(H_\omega)=2$. We have
\begin{align*}
& H_\omega\sigma_3\Phi_\omega=0,\quad H_\omega\partial_\omega\Phi_\omega
=-\Phi_\omega.
\end{align*}
Since $H_\omega^*=\sigma_3H_\omega\sigma_3$, we have
$N_g(H_\omega^*)=\operatorname{span}\{\Phi_\omega,
\sigma_3\partial_\omega\Phi_\omega\}$. Let $\xi(\omega)$ be a  real
eigenfunction with eigenvalue   $\lambda(\omega)$. Then we have
$$
H_\omega\xi(\omega)=\lambda(\omega)\xi(\omega),\quad
H_\omega\sigma_1\xi(\omega)=-\lambda(\omega)\sigma_1\xi(\omega).$$
Notice that $\langle  \xi,\sigma _3 \xi\rangle >0$ since
$\langle  \sigma_3 H_\omega \cdot,\cdot \rangle $ is positive definite on
$  N_g^\perp (H_\omega^*)$.
For $\omega\in\mathcal{O}$, we have the $H_\omega$-invariant Jordan
block decomposition
\begin{align}  \label{eq:spectraldecomp} &
L^2(\R ,\C^2)=N_g(H_\omega)\oplus \big (\oplus_{\pm}N(H_\omega\mp
\lambda(\omega)) \big)\oplus L_c^2(H_\omega),
\end{align}
where $L_c^2(H_\omega):=
\left\{N_g(H_\omega^\ast)\oplus(\oplus _{\pm  }N(H_\omega^\ast \mp
\lambda(\omega))\right\} ^\perp .$
Correspondingly, we set
\begin{align}
  \label{eq:decomp2}
& R(t) =z(t)\xi(\omega(t))+\overline{z}(t) \sigma_1\xi(\omega(t))+f(t),\\
\label{eq:decomp3}
& R(t)\in {} N_g^\perp(H_{\omega(t)}^*)\quad\text{and}\quad
f(t)\in L_c^2(H_{\omega(t)}).
\end{align}
There is a Taylor expansion at $R=0$ of the
nonlinearity $O(R^2)$ in \eqref{eqR} with $R_{m,n}(\omega  ,x) $ and
$A_{m,n}(\omega ,x ) $ real vectors  and matrices rapidly decreasing
in $x$: $ O(R^2)=$
\begin{equation*}
 \sum _{ 2\le  m+n \le 2N+1} R_{m,n}(\omega ) z^m  \bar z^n+
\sum _{1\le  m + n \le  N} z^m  \bar z^n A_{m,n}(\omega ) f+
O(f^2+|z|^{2N+2}) .\end{equation*}
Then,
\begin{equation}\label{eq:Rsplit}\begin{aligned}& \im f_t=\left ( H _{\omega (t)}+\sigma _3 \dot \gamma \right )f + \sigma _3
\dot \gamma \Phi (\omega )- \im  \dot \omega \partial _\omega \Phi (t)
+
  (z \lambda  (\omega ) -\im \dot z ) \xi  (\omega ) \\ &
- (\bar z \lambda  (\omega )+i\dot {\bar z }) \sigma _1\xi (\omega )
  +\sigma _3 \dot \gamma (z  \xi + \bar z  \sigma _1 \xi )
-i \dot \omega (z  \partial _\omega \xi + \bar z  \sigma _1
\partial _\omega \xi ) \\ & + \sum _{ 2\le  m+n \le 2N+1} z^m \bar
z^nR_{m,n}(\omega ) + \sum _{1\le  m + n \le N} z^m \bar z^n
A_{m,n}(\omega ) f+\\ & + O(f^2)+ O_{loc}(|z  ^{2N +2}|) ,
\end{aligned}
\end{equation}
where by $O_{loc}$ we mean that the there is a factor $\chi (x)$
rapidly decaying   to 0   as $|x|\to \infty $. Taking inner products
of the equation with the generators of $N_g(H_\omega ^\ast )$ and
of $\ker (H_\omega ^\ast -\lambda  )$, we obtain modulation and
discrete modes equations ($q'(\omega ):=\frac{d\| \phi _\omega \|
_2^2}{d\omega  } $):
\begin{equation}\label{eq:discrete}\begin{aligned}&
\im \dot \omega q'(\omega ) =\langle \mathcal{X}, \Phi \rangle \, ,
\,  \dot \gamma q'(\omega ) =\langle \mathcal{X}, \sigma _3
\partial _\omega \Phi \rangle  \, , \,  \im \dot z -\lambda (\omega )z =\langle
\mathcal{X}, \sigma _3 \xi  \rangle \, , \, \\ & \mathcal{X}:= \sigma
_3 \dot \gamma (z \xi + \bar z  \sigma _1 \xi ) - \im \dot \omega (z
\partial _\omega \xi + \bar z  \sigma _1
\partial _\omega \xi )  +   \sum _{    m+n =2}^{ 2N+1}
 z^m \bar z^nR_{m,n}(\omega )\\ & +  \big ( \sigma _3 \dot \gamma + \im \dot \omega
\partial _\omega P_c+ \sum _{  m + n =1}^{ N} z^m \bar z^n
A_{m,n}(\omega )  \big ) f   + O(f^2)+ O_{loc}(|z ^{2N +2}|) .
\end{aligned}
\end{equation}

We now go through the dispersive estimates. The proofs are in \cite{Cu5}.
We call admissible a pair   $(p,q)$  s.t.
\begin{equation}\label{admissiblepair}  2/p+1/q= 1/2\,
 , \quad   p\ge 4\, , \quad
q\ge 2.
\end{equation}

\begin{theorem}[Strichartz estimates]
\label{Thm:Strichartz} For   $k =[0,2]$ there exist  positive
numbers $ C(\omega , k , p)$  and $  C(\omega , k , p_1, p_2)$  upper semicontinuous in their arguments such that:
  {\item {(a)}}
 for any $f\in
L^2_c( {\omega })$ and any admissible pair $(p,q)$ with $p>4$ we have
\begin{equation}\label{stri1}\|e^{-\im tH_{\omega }} f\|_{L_t^pW_x^{k,q } }\le C \|f\|_{H^k};\end{equation}
 {\item {(b)}}
  for any $g(t,x)\in
S(\R^2)$ and any   two admissible pairs $(p_j,q_j)$  for $j=1,2$    with $p_j>4$ we have
\begin{equation}\label{stri2}
\|\int_{0}^te^{-\im (t-s)H_{\omega }} P_c( {\omega
})g(s,\cdot)ds\|_{L_t ^{p_1}W_x^{k,q_1 }} \le
C \|g\|_{L_t^{p_2'}W_x^{k,q_2' } }.
\end{equation}
In the case $k=0$, we can include also case $p=4$ in \eqref{stri1}
and $p_j=4$ for any of $j=1,2$ in \eqref{stri2}.
\end{theorem}
For the  proof see \cite{Cu5}. The case
 $k>0$ requires interpolation. The case $k=0$ is like the one for
 $e^{-\im t\partial _x^2}$.  Specifically, we can use   dispersive
 estimates, see \cite{KS}, and an appropriate version of the so called
 $TT^*$ argument.  In particular, this yields the  $L_t^4L_x^{ \infty
 }$  bound, which is not reached in \cite{KS}.  See \cite{DMW} for how
 to extend such results to Schr\"odinger operators, $H$, formed by
 singular perturbations of the Laplacian with $k \leq 1$.

\begin{lemma}
  \label{lem:resolvent} Fix $\tau >3/2$.

{\item {(1)}} There exists $C=C(\tau ,\omega )$, upper
semicontinuous in $\omega $ such that  for any $\varepsilon \neq 0$,
$$\| R_{H_\omega }(\lambda +\im \varepsilon )P_c(H_\omega )
u\| _{L^2_\lambda L^{2,-\tau }_x}\le  C \| u\| _{L^2}. $$

{\item{(2)}} For any $u\in L^{2, \tau }_x $ the following limits exist:
$$ \lim _{\epsilon \searrow 0}R_{H_\omega }(\lambda \pm \im \varepsilon )
u= R_{H_\omega }^\pm  (\lambda ) u  \text{ in $C^0(\sigma
_e(H_\omega ),L^{2, -\tau }_x)$}.$$

{\item {(3)}} There exists $C=C(\tau ,\omega )$, upper
semicontinuous in $\omega $ such that
$$
  \|   R_{H_\omega }^\pm  (\lambda
)P_c(H_\omega )   \| _{B( L^{2,\tau }_x, L^{2,-\tau }_x)} < C
\langle \lambda \rangle ^{-\frac{1}{2}} .$$

{\item {(4)}} Given any
$u\in L^{2, \tau }_x $ we have
$$P_c(H_\omega )u=\frac{1}{2\pi \im }\int _{\sigma _e(H_\omega )}
(R_{H_\omega }^{+}(\lambda  )-R_{H_\omega }^{-}(\lambda  ))  u\,
d\lambda .$$
\end{lemma}

These are consequences of the fact that $\sigma _e(H_\omega )$ does
not contain eigenvalues and that $\pm \omega $ are not resonances,
and of the theory on plane waves and representation of the resolvent
in \cite{KS}. In fact, of a much simpler version than \cite{KS}, due
to the fact that $H_\omega $ is a small perturbation of $\sigma _3
(-\partial _x^2+V +\omega ) $.

Claim (a) of the  following smoothing lemma is a consequence of
Lemma \ref{lem:resolvent} by \cite{kato}, while   (b) follows from
(a) by duality.

 \begin{lemma}
  \label{lem:smoothing}  For any $k, \tau >3/2$,  $\exists$
  $C=C(\tau ,k,\omega )$ upper
semicontinuous in $\omega $ such that:
   {\item {(a)}}
  for any $f\in S(\R )$,
$$
\| e^{-\im tH_{\omega }}P_c(H_\omega )f\| _{L_{  t}^2 H_x^{k, -\tau}}
\le
 C\|f\|_{H ^{k}} .
  $$
  {\item {(b)}}
  for any $g(t,x)\in
 {S}(\R ^2)$
$$ \left\|\int_\R e^{\im tH_{\omega }}
P_c(H_\omega )g(t,\cdot)dt\right\|_{H^k_x} \le C\| g\|_{L_{  t}^2
H_x^{k, \tau}}.
$$
\end{lemma}

\begin{lemma}
  \label{lem:endpoint} For any $k, \tau >3/2$,  $\exists$
  $C=C(\tau ,k,\omega )$ as
above such that $\forall$ $g(t,x)\in {S}(\R^2)$
$$  \left\|  \int_0^t e^{-\im (t-s)H_{\omega
}}P_c(H_\omega )g(s,\cdot)ds\right\|_{L_{  t}^2 H_x^{k, -\tau}} \le
C\|  g\|_{L_{  t}^2 H_x^{k, \tau}}.
$$
\end{lemma}
\proof   To get this proof there is no need of Lemma 11 \cite{M1}
or of the analogous result, Lemma 2 in \S 7, in  \cite{PelinovskyStefanov}.
 We just use Plancherel and H\"older inequalities and   (3)
Lemma
 \ref{lem:resolvent}:
\begin{equation*} \begin{aligned} &
\| \int _{0}^t e^{-\im (t-s)H_{\omega }}P_c(H_\omega )g(s,\cdot)ds\|_{
L_{t }^2L_{ x}^{2,-\tau }} \le \\& \le  \| R_{H_\omega }^+(\lambda
)P_c (H_\omega )
  \widehat{ \chi }_{[0,+\infty )}\ast _\lambda
   \widehat{ g}(\lambda,x)\|_{L_{t }^2L_{ x}^{2,-\tau } }   \le \\& \le
   \left\| \,
\|   R_{H_\omega }^+ (\lambda )P_c  (H_\omega ) \| _{B(
L^{2,\tau}_x, L^{2,- \tau}_x)} \|
   \widehat{ \chi }_{[0,+\infty )}
   \ast _{\lambda } \widehat{g} (\lambda,x) \|_{L_{ x}^{2, \tau }}\, \right\|_{L^2_\lambda}
\\ \le &
  \|  R_{H_\omega }^+ (\lambda
)P_c(H_\omega )   \| _{L^\infty _\lambda (\Bbb R ,B( L^{2,\tau}_x,
L^{2,-\tau}_x))}\| g\|_{L_{t }^2L_{ x}^{2, \tau } }  \le C \|
g\|_{L_{t }^2L_{ x}^{2, \tau } } .
\end{aligned}
\end{equation*}
\qed

\begin{lemma}
  \label{lem:L4} $k$ and  $\tau >3/2$  $\exists$
  $C=C(\tau ,k,\omega )$ as
above such that $\forall$ $g(t,x)\in {S}(\R^2)$
$$
\left\|\int_0^t e^{-\im (t-s)H_{\omega }}P_c(H_\omega )g(s,\cdot)ds
  \right\|_{
L_t^\infty L_x^2\cap  L^{4 }_{t}(\R ,W ^{k,\infty }_x) } \le
C\|g\|_{L_t^2H_x^{k,\tau}}.
 $$
\end{lemma}
\proof For $g(t,x)\in S(\R^2)$ set
$$T g(t)=\int _0^{+\infty} e^{-\im (t-s)H_{\omega }}
P_c(H_\omega )g(s) ds.$$ Lemma \ref{lem:smoothing} (b) implies $f:=\int _0^{+\infty}
e^{\im sH_\omega }P_c( \omega )g(s)ds\in L^2(\R)$. Then Lemma \ref{lem:L4}
 is a direct consequence of \cite{CK}.\qed

 The following lemma can be proved in a way similar to Lemma B.1 \cite{Cu5}.
 \begin{lemma}
  \label{lem:trick}  The following operators  $P_\pm (\omega )$ are
well-defined:
\begin{equation*}   \begin{aligned} & P_+(\omega )u =\lim _{\epsilon \to 0^+}
 \frac 1{2\pi \im }
\lim _{M \to +\infty} \int _\omega ^M \left [ R_{H_\omega }(\lambda +\im \epsilon
)- R_{H_\omega }(\lambda -\im \epsilon )
 \right ] ud\lambda , \\ &
P_-(\omega )u =\lim _{\epsilon \to 0^+}
 \frac 1{2\pi \im }
\lim _{M \to +\infty} \int _{-M }^{-\omega } \left [ R_{H_\omega }(\lambda
+\im \epsilon )- R_{H_\omega }(\lambda -\im \epsilon ) \right ] ud\lambda .
\end{aligned}
\end{equation*}
  For any
$M>0$ and $N>0$ and for $C=C (N,M,\omega  )$ upper semicontinuous in
$\omega >\omega ^*$, we have $$  \|     (P_+(\omega
 )-P_-(\omega  )-P_c(\omega  )\sigma _3) f\|  _{L^{2,M} }\le
C  \|      f\|  _{L^{2,-N} }. $$
\end{lemma}

 \section{Normal form expansion}
\label{sec:normal}

Here we repeat the theory in \cite{cuccagnamizumachi,Cu5} which is
somewhat more elementary than \cite{Cu1}, but still adequate in our
setting.

We consider $N\in \mathbb{N}$ such that for any $t\ge 0$, for $\rho
_1(t)$ and for the corresponding $\omega (t):=\omega (\rho _1(t))$
and for $\lambda (t)= \lambda (\omega (\rho _1(t)))$ we have
\begin{equation}\label{eq:Resonance}N \lambda (t)<\omega
(t)<(N+1) \lambda (t).\end{equation} Notice that $\frac{\lambda
(\omega (\rho _1 ))}{\omega (\rho _1 )}$ is for $\rho _1 \ge 0$ a
continuous and strictly increasing function, equal to 0 at $\rho _1=0$. This means that for
$\rho _1 > 0$ small, it has no values in $\mathbb{N}$. By continuity
and orbital stability, we can then assume  \eqref{eq:Resonance}
for all $t$.

 \subsection{Changes of variables on $f$}
\label{subsec:normal f}

For the $N$ of \eqref{eq:Resonance}
 we   consider  $k=1,2,...N$ and
set $f=f_k$ for $k=1$. The other  $f_k$ are defined below. In the
ODE's  there will be error terms    of the form
$$E_{ODE}(k)=   O( |z |^{2 N+2 }
 )+O(  z^{N+1}   f _{k}  )  +O(f^2_{k})+
O(|f _{k}|^5 ).$$ In the PDE's there will be error
terms   of the form
$$E_{PDE}(k)=   O_{loc}( |z |^{N+2}
 )  +O_{loc}(  z f _{k}  )+O_{loc}(f^2_{k})+
O( |f _{k}|^4f _{k} ).$$ In the right hand sides of the
equations  (2.3-4) we substitute $\dot \gamma $ and $\dot \omega $
using the modulation equations. We repeat the procedure a sufficient
number of times until we can write for $k=1$,  $ f_1=f$ and $q'(\omega )= \frac{d\| \phi _\omega \| _2^2}{d\omega  },$
\begin{equation}\label{eq:stepk} \begin{aligned} &   \im \dot \omega q'(\omega )
= \langle   \sum _{   m +n =2}^{ 2N+1 } z^m \bar z^n\Lambda
_{m,n}^{(k)}(\omega ) +   \sum _{   m +n =1}^{ N } z^m \bar z^n
A_{m,n}^{ (k) }(\omega ) f_{k}  \\ & +E_{ODE}(k) , \Phi (\omega ) \rangle
  \\  & \im \dot z -\lambda  z = \langle
\text{ same as above }, \sigma _3 \xi (\omega )  \rangle
\\ & \im \partial _t f_{k}= \left ( H _{\omega}  +\sigma _3 \dot \gamma
\right  )f_{k} + E_{PDE}(k)+ \sum _{k+1\le  m +n \le N +1}z^m \bar
z^n R_{m,n}^{(k)}(\omega ),
\end{aligned}\end{equation}
 with  $A_{m,n}^{(k)} $, $R_{m,n}^{(k)}     $  and $\Lambda
_{m,n}^{(k)} (\omega , x) $ real
 exponentially  decreasing  to 0 for $|x|\to \infty$ and continuous in $(\omega , x) $.
Exploiting $|(m-n)\lambda (\omega )|<\omega $ for $  m+n \le N$,
$m\ge 0$, $n\ge 0$,
 we define inductively $f_k$ with $k\le N$ by
\begin{eqnarray*}
f_{k } & = & \sum _{  m+n=k}z^m \bar z^n \Psi _{m,n}(\omega )
+f_{k -1}  \text{ for }\Psi _{m,n}(\omega ) \\
& := &R_{H_{\omega}}((m-n)\lambda (\omega ) ) R_{m,n}^{(k-1)}.
\end{eqnarray*}
Notice that if $ R_{m,n}^{(k-1)}   (\omega
 ,x)$ is real
 exponentially  decreasing  to 0 for $|x|\to \infty$, the same is
 true for $\Psi _{m,n}(\omega )$ by $|(m-n)\lambda (\omega )|<\omega $. By induction $f_k$ solves the above equation with
the above notifications.

We are now ready to state the result which directly implies
Theorem \ref{thm:as stab}.
\begin{theorem}\label{thm:dip}
 Assume (H1)--(H5).
Let $u$ be a solution of \eqref{NLS}, $U={}^t\!(u,\overline{u})$,
and let $\Psi _{m,n}(\omega)$ be as above. Then if
$\epsilon_0$ in Theorem \ref{thm:as stab} is sufficiently small,
there exist $C^1$-functions $\omega(t)$ and $\theta(t)$, a constant
$\omega_+>\omega ^*$ such that we have
$\sup_{t\ge0}|\omega(t)-\omega_0|=O(\epsilon)$, $\lim _{t\to
+\infty  } \omega(t)=\omega_+$ and we can write
\begin{align*}
U(t,x) = & e^{\im  \theta (t)\sigma_3}
\left(\Phi_{\omega (t)}(x)+z(t)\xi(\omega(t))
+\overline{z(t)}\sigma_1\xi(\omega(t))\right)
\\ & +
e^{i \theta (t)\sigma_3}
\sum_{ 2\le m+n\le N }\Psi_{m,n}(\omega(t))
z(t)^m\overline{z(t)}^n+e^{i \theta (t)\sigma_3}f_N(t,x),
\end{align*}
\begin{gather*}\text{with }
\| z(t)\|_{ L_t ^{2N+2}}^{N+1} +
\|f_N(t,x)\|_{L^\infty  _t H ^{1 }_x  \cap L^5 _t W ^{1,10}_x
 \cap L^4_t L^\infty _x }
\le C \epsilon.
\end{gather*}
Furthermore, there exists $f_+\in H^1(\R ,\C^2)$ such that
\begin{equation}  \label{eq:scatt}  \lim_{t\to +\infty} \left \|  e^{\im \theta (t) \sigma _3}f_N(t) -
 e^{\im t \partial ^2_x   \sigma_3}{f}_+  \right \|_{H^1}=0.\end{equation}
\end{theorem}
Obviously the scattering result \eqref{eq:scatt} holds also for $t\to -\infty$. We do not prove this lemma explicitly, but we recall
Lemma 4.3 in \cite{Cu5} which states:

\begin{lemma}
  \label{lem:cont} There are fixed  constants $C_0$ and $C_1$ and
$\epsilon _0>0$  such that   for any $0<\epsilon \le \epsilon _0$
  if we have
\begin{equation}\label{eq:cont1}\|  { z} \| _{L^{2N+2 }_t}^{N+1} \le 2C_0\epsilon \quad \&
 \quad \| f_N \| _{L^\infty  _tH ^{1 }_x \cap L^5 _tW ^{1,10}_x
 \cap L^4_tL^\infty _x\cap L^2_tH^{1,-2}_x } \le 2C_1\epsilon,
  \end{equation}
then we obtain the improved inequalities
\begin{align}\label{eq:cont2}
 &
\| f_N \| _{L^\infty  _tH ^{1 }_x \cap L^5 _tW ^{1,10}_x \cap
L^4_tL^\infty _x \cap  L^2_tH^{1,-2}_x} \le C_1\epsilon  ,
\\  \label{eq:cont3}  & \|  {z} \| _{L^{2N+2 }_t}^{N+1} \le C_0\epsilon .
\end{align}
\end{lemma}
We sketch only the main steps of the proof. First of all
 we rewrite  the   equation for $f_N$.
Set $\omega (0):=\omega (\rho _1(0))$
 and write

\begin{equation*} \begin{aligned} & \im   \partial _t
f_{N}-\left ( H_{\omega (0)}  +  \sigma _3 \left  (   \dot \gamma
+\omega -\omega (0)  \right  )\right )   f_{N}
    =\sum _{  m+n= N+1} z^m \bar
z^n   R_{m,n}^{(N)} (\omega   )\\& + \widetilde{E}_{PDE}(N) \text{
where } \widetilde{E}_{PDE}(N):= E_{PDE}(N)  +  (  \mathcal{V}
_\omega
 -\mathcal{V} _{\omega (0)}) f_{N}.
\end{aligned}\end{equation*}
It is easy to see that \eqref{eq:cont2} for $f_N $, which is the same of $P_c(\omega  )f_N$,  or for $P_c(\omega (0))f_N$, are equivalent. This because  $P_c(\omega  )-P_c(\omega (0)) =P_d(\omega (0)) -P_d(\omega  ) $ is a small and smoothing operator. We then write for $\varphi (t)=\varphi = \gamma
+\omega -\omega (0) $

\begin{equation*} \begin{aligned} & \im   \partial _t
P_c(\omega (0))   f_{N}-\left ( H_{\omega (0)}  +    \varphi  \left (P_+(\omega (0))  - P_-(\omega (0)) \right  )\right )   P_c(\omega (0)) f_{N}
    =\\& \sum _{  m+n= N+1} z^m \bar
z^n   P_c(\omega (0))R_{m,n}^{(N)} (\omega   )  + \\& P_c(\omega (0))\widetilde{E}_{PDE}(N) + \varphi \left (   P_+(\omega (0))  - P_-(\omega (0)  ) -P_c(\omega (0))\sigma _3     \right )    f_{N} .
\end{aligned}\end{equation*}
It turns out that the term in the second line is the main one of the rhs
and that the norms of $f_N$ in \eqref{eq:cont2} can be controlled by
$\| f _N (0) \| _{H^1}+\|  {z} \| _{L^{2N+2 }_t}^{N+1}$. In particular for
$\widetilde{E}_{PDE}(N)$ we refer to Lemma 4.5 \cite{Cu5}. The very last term
in the equation is controlled using the fact that $\varphi $ is small, Lemmas
\ref{lem:endpoint}, \ref{lem:L4} and  \ref{lem:trick}.

\subsection{A further change of variable in $f$}
\label{subsec:normal g}

In the argument there is need for a decomposition of $f_N$, namely setting
\begin{equation}\label{eq:g}
 f_{N}=-\sum _{  m+n=N+1} z^m \bar z^n   R_{  H_{\omega (0)}      }^{+}((m-n)
\lambda     )
 P_c(\omega (0))R_{m,n}^{(N)}   (\omega  )  +  g ,
\end{equation}
 where if $| \Lambda  |< \omega (0)$, we set $ R_{  H_{\omega (0)}      }^{+}(\Lambda    ) =R_{  H_{\omega (0)}      } (\Lambda    )$.
 The function $g$ satisfies an equation of the form
 \begin{equation*} \begin{aligned} & \im \partial _t P_c(\omega (0)) g=\left ( H_{\omega (0)}    + \varphi
(P_+(\omega (0)) -P_-(\omega(0)) )\right ) P_c(\omega (0)) g +
\\& +  \sum _{\pm } O(\epsilon |z|^{N+1}
)R_{H_{\omega (0)}} ^+(\pm  (N+1) \lambda (\omega (0))  ) R_\pm
  +P_c(\omega (0))\widehat{ {E}}_{PDE}(N)
  \end{aligned}\end{equation*}   $R_+ :=R_{ N+1,0}^{(N)} $,  $R_- :=R_{
0, N+1 }^{(N)} $ and $\widehat{ {E}}_{PDE}(N):=\widetilde{
{E}}_{PDE}(N)+ O_{loc}(\epsilon z^{N+1})$.

\begin{lemma}
  \label{lem:estimate g} Assume the hypotheses of Lemma \ref{lem:cont}.
Then, there exists a fixed    $C_0 =C(\omega (0))$   such that for a fixed $S$ sufficiently large
\begin{equation} \label{bound:auxiliary}\| g
\| _{L^2_tL^{2,-S}_x}\le C_0 \epsilon + O(\epsilon
^2).\end{equation}
\end{lemma}
See Lemma 4.6 \cite{Cu5}.

 \subsection{Change of $\omega $ and $z$.}
\label{subsec:normal z}

 Consider now the equations of $\omega $ and $z$ in \eqref{eq:stepk}.
 Then, we have the following

 \begin{lemma}\label{lemma:change z}  There is a change of variables
 \begin{equation}\label{eq:change z} \begin{aligned} &\widetilde{\omega }=\omega +q(\omega , z,\bar z) +\sum _{1\le  m
+n \le N }z^m  \bar z^n\langle f_N,A _{mn}(\omega )\rangle ,\\&
\zeta =z  +p (\omega , z,\bar z) +\sum _{1\le  m +n \le N
}z^m \bar z^n\langle f_N,B _{mn}(\omega )\rangle ,
\end{aligned}\end{equation} with $p (\omega , z,\bar z)=\sum p_{ m,n}(\omega )
z^m\bar z^n$ and $q(z,\bar z)=\sum q_{ m,n}(\omega ) z^m\bar z^n$
polynomials in $(z,\bar z)$ with real coefficients and $O(|z|^2)$
 near 0, such that we get for $   a_{ m  }(\omega )  $ real

\begin{equation} \label{finalsyst}\begin{aligned} & \im  \dot {\widetilde{\omega}} =  \langle  {E}_{PDE}(N) , \Phi \rangle   \\ &
 \im \dot {\zeta }-\lambda  (\omega )\zeta  =
 \sum _{ 1\le  m \le N} a_{ m  }(\omega )|\zeta  | ^{2m }
 \zeta +\langle  E_{ODE}(N) , \sigma _3 \xi   \rangle \\
& +
  \overline{\zeta}^N \langle  A_{0,N}^{(N)}(\omega ) f _{N} , \sigma _3
\xi   \rangle  .
\end{aligned}\end{equation}
\end{lemma}
\proof The proof is elementary and goes as follows, see \cite{cuccagnamizumachi} for
details. We consider recursively for $\ell=0,..., 2N$ with $z_0=z$ equations

\begin{equation}\label{zell}\begin{aligned} &
 \im \dot  z_\ell-\lambda z_\ell= \sum _{ 1\le l \le \ell } a_{ l,\ell  }(\omega )|z_\ell| ^{2 l}z_\ell+ \sum _{ \ell +2\le m+n \le 2N+1} z^m _\ell \overline{z}^n_\ell
 \alpha ^{(\ell )}_{m,n}(\omega ) \\& + \sum _{ \ell +1\le m+n \le  N } z^m _\ell \overline{z}^n_\ell
 \langle A ^{(\ell )}_{m,n}(\omega ) , f_N \rangle +E_{ODE}(\ell) \end{aligned}
\end{equation}
with  $\alpha ^{(\ell )}_{m,n}(\omega )\in \R$ and
$A ^{(\ell )}_{m,n}(\omega )\in \mathcal{S}(\R ^3, \R^2)$. Suppose this holds for $\ell < 2N$. Then set
\begin{equation*} \begin{aligned} &
   z_  {\ell +1}=z_  {\ell } + \sum _{   m+n = \ell +2} z^m _\ell \overline{z}^n_\ell
 \beta ^{(\ell )}_{m,n}(\omega )   + \sum _{   m+n =\ell +1 } z^m _\ell \overline{z}^n_\ell
 \langle B ^{(\ell )}_{m,n}(\omega ) , f_N \rangle ,  \\& \beta ^{(\ell )}_{m,n}(\omega ):=  \frac{ \alpha ^{(\ell )}_{m,n}(\omega )}{(m-1-n) \lambda} \text{ for $m\neq n+1$, } \beta ^{(\ell )}_{n+1,n}(\omega )=0 ,\\&  B ^{(\ell )}_{m,n}(\omega ) = -R _{H ^{*}_\omega }  ((m-1 -n) \lambda) A ^{(\ell )}_{m,n}(\omega ) \text{ for $\ell<N $, } \\&  B ^{(N )}_{m,n}(\omega ) = -R _{H ^{*}_\omega }  ((m-1 -n) \lambda) A ^{(\ell )}_{m,n}(\omega ) \text{ for $(m,n)\neq (0,N) $,} \end{aligned}
\end{equation*}
where $B ^{(N )}_{0,N}(\omega ) =0 $ otherwise and $B ^{(\ell )}_{m,n}(\omega ) =0$
for $\ell >N$. This yields for $ \zeta =z _{2N} $ the desired result.

 We substitute in the  equation for $\omega $ in \eqref{eq:stepk}  (for $k=N$) $z$ with $\zeta $ inverting the second equation in \eqref{eq:change z}.
  We consider recursively for $\ell=0,..., 2N+1$ with $\Omega _0=\omega$ equations
\begin{equation*} \begin{aligned}
 \im \dot  \Omega _\ell &=    \sum _{ \ell +2\le m+n \le 2N+1} z^m _\ell \overline{z}^n_\ell
 \gamma ^{(\ell )}_{m,n}(\omega )  + \sum _{ \ell +1\le m+n \le  N } z^m _\ell \overline{z}^n_\ell
 \langle \Gamma ^{(\ell )}_{m,n}(\omega ) , f_N \rangle \\ &
 +E_{ODE}(\ell)
\end{aligned}
\end{equation*}
with  $\gamma ^{(\ell )}_{m,n}(\omega )\in \R$ and
$\Gamma ^{(\ell )}_{m,n}(\omega )\in \mathcal{S}(\R ^3, \R^2)$.
 Suppose this holds for $\ell < 2N+1$. Then set
\begin{equation*} \begin{aligned} &
  \Omega _  {\ell +1}=\Omega _  {\ell } + \sum _{   m+n = \ell +2} z^m _\ell \overline{z}^n_\ell
 \delta  ^{(\ell )}_{m,n}(\omega )   + \sum _{   m+n =\ell +1 } z^m _\ell \overline{z}^n_\ell
 \langle \Delta  ^{(\ell )}_{m,n}(\omega ) , f_N \rangle ,  \\& \delta ^{(\ell )}_{m,n}(\omega ):=  \frac{ \gamma ^{(\ell )}_{m,n}(\omega )}{(m -n) \lambda} \text{ for $m\neq n $, } \delta ^{(\ell )}_{n ,n}(\omega )=0 \\&  \Delta ^{(\ell )}_{m,n}(\omega ) = -R _{H ^{*}_\omega }  ((m -n) \lambda) \Gamma ^{(\ell )}_{m,n}(\omega ) \text{ for $\ell \le N  $, }   \end{aligned}
\end{equation*}
  $\Delta^{(\ell )}_{m,n}(\omega ) =0$
for $\ell >N$. This yields for $  \widetilde{\omega} =\Omega  _{2N+1} $ the desired result.
 \qed

 By  \eqref{eq:cont1} we have $\| \dot {\widetilde{\omega}}\| _{L^1_t}=O(\epsilon ^2)$.

\begin{remark}
\label{rem:approxiamtion} Setting $\widetilde{\omega} (t)\equiv \widetilde{\omega} (0)$,      $f_N\equiv 0$  and considering the equation
\begin{equation*} \label{finalsyst1}\begin{aligned} &
 \im \dot {\zeta }-\lambda  (\omega )\zeta  =
 \sum _{ 1\le  m \le N} a_{ m  }(\omega )|\zeta  | ^{2m }
 \zeta
\end{aligned}\end{equation*}  yields a finite dimensional approximation
of the NLS. We do not check here the time span when the solutions
of this approximation  are good approximations of solutions of the full NLS.
Nonetheless we recall that in the literature, see for example \cite{BP2,GS}, are displayed solutions of \eqref{finalsyst} s.t.
approximately $|\zeta (t)|\approx \dfrac{|\zeta (0)|}{ (|\zeta(0)|^{2N} \, N\, \Gamma \, t+
1)^{\frac{1}{2N}}}$, with this approximation valid for  $t\ll |\zeta(0)|^{-2N} $.

\end{remark}

\section{The Fermi golden rule}
\label{sec:FGR}

Substituting in the equation for $\zeta$ the variable  $f_N$ with \eqref{eq:g} to get
\begin{equation*} \begin{aligned} &
 \im  \dot {\zeta}-\lambda  (\omega )\zeta  =
 \sum _{ 1\le m\le N}  {a}_{ m  }(\omega )|\zeta | ^{2m }
 \zeta +\langle  E_{ODE}(N) , \sigma _3 \xi   \rangle   -\\&
 -
  |\zeta  |^{2N}\zeta  \langle   {A}_{0,N}^{(N)}(\omega )
  R_{H_{\omega (0)}}^+((N+1) \lambda  (\omega (0)) )P_c(\omega
_0)R_{N+1 ,0}^{(N)}(\omega  )     ,\sigma _3\xi   \rangle
\\& +  \overline{\zeta} ^N \langle  {A}_{0,N}^{(N)}(\omega ) g , \sigma _3 \xi   \rangle
\end{aligned}\end{equation*}
with $ {a}_{ m  }$, $ A _{0,N}^{(N)}$ and
$R_{N+1,0}^{(N)}$ real. Set
\begin{equation*} \begin{aligned} &\Gamma  (\omega ,
\omega (0))= \Im \left (\langle  {A}_{0,N}^{(N)}(\omega )
  R_{H_{\omega  (0)}} ^{+}((N+1) \lambda (\omega (0)) \right. \\
& \left. \times P_c(\omega
(0))R_{N+1,0}^{(N)}(\omega  ) \sigma _3 \xi   (\omega )\rangle \right
)\\& =\pi
  \langle   {A}_{0,N}^{(N)}(\omega )
  \delta (H_{\omega (0)}- (N+1) \lambda (\omega (0))  )P_c(\omega (0))R_{N+1,0}^{(N)}(\omega  ) \sigma _3 \xi   (\omega )\rangle .
 \end{aligned}\end{equation*}
Now we assume the following: \begin{itemize}
\item[(H5)] There is a
fixed constant $\Gamma >0$ such that $|\Gamma  (\omega , \omega
)|>\Gamma .$

\end{itemize}
It is then easy to see, for example Corollary $4.7$ \cite{Cu5}, that in fact $ \Gamma
(\omega , \omega ) >\Gamma  $, but we will not use this here. By
continuity,
we can assume $|\Gamma (\omega ,\omega (0) )|>\Gamma /2.$ Then, we
write
\begin{equation*} \begin{aligned} &  \frac{d}{dt} \frac{|\zeta |^2}{2} =    -\Gamma  (\omega , \omega  (0))
  |\zeta |  ^{2N+2} +\Im \left (
\langle  {A}_{0,N}^{(N)}(\omega )f _{N+1} , \sigma _3 \xi
(\omega )\rangle \overline{\zeta} ^{N+1}   \right ) \\&+ \Im
\left ( \langle E_{ODE}(N) , \sigma _3 \xi  (\omega ) \rangle
\overline{\zeta} \right ) .
\end{aligned}\end{equation*}
For $A_0$   an
upper bound  of the    constants $A_0(\omega )$ of Theorem \ref{thm:orbstab},
we get
$$\frac{\Gamma}{2} \| \zeta\| _{L^{2N+2}_t}^{2N+2}\le  A_0\epsilon ^2+2c(\omega (0))\epsilon
\| \zeta\| _{L^{2N+2}_t}^{N+1} + o(\epsilon ^2).$$ Then we can
pick    $C_0= 2 ( A_0+2c(\omega (0) +1))/\Gamma $ in Lemma \ref{lem:cont}
and this proves that \eqref{eq:cont1}
implies \eqref{eq:cont3}. Furthermore $\zeta(t)\to 0$ since
$\frac{d}{dt}\zeta (t)=O(\epsilon ),$ again see \cite{Cu5} for more in
depth discussion.

\appendix

\section{Finite Dimensional Dynamics}
\label{sec:fd}

We plug the ansatz
\begin{eqnarray*}
u(x,t) = c_0 (t) \psi_0 + c_1 \psi_1 + R(x,t)
\end{eqnarray*}
into \eqref{NLS}, where
\begin{eqnarray*}
H \psi_j & = & (-\p_x^2 + V) \psi_j \\
& = & - \omega_j \psi_j
\end{eqnarray*}
for $j = 0,1$.  As a result, we have the equation
\begin{eqnarray*}
i \dot{c}_0 \psi_0 + i \dot{c}_1 \psi_1 + i R_t (x,t) & = & -\omega_0
c_0 \psi_0 - \omega_1 c_1 \psi_1 + H R \\
&& - |c_0 \psi_0 + c_1 \psi_1 +
R|^4 (c_0 \psi_0 + c_1 \psi_1 + R).
\end{eqnarray*}
The nonlinear contribution is then given by
\begin{eqnarray*}
( c_0 \psi_0 \!\!\!\!\! &+& \!\!\!\!\!  c_1 \psi_1)^3 (\bar{c}_0 \psi_0 + \bar{c}_1 \psi_1)^2 +
\mathcal{O} (R) = \\
 (c_0^3 \psi_0^3 \!\!\!\!\! &+& \!\!\!\!\!  3 c_0^2 c_1 \psi_0^2 \psi_1 + 3 c_0 c_1^2 \psi_0
\psi_1^2 + c_1^3 \psi_1^3) (\bar{c}_0^2 \psi_0^2 + 2 \bar{c}_0
\bar{c}_1 \psi_0 \psi_1 + \bar{c}_1^2 \psi_1^2) \\
&& + \mathcal{O} (R)  = \\
(c_0^3 \bar{c}_0^2) \psi_0^5 \!\!\!\!\! &+& \!\!\!\!\!  (3 c_0^2 \bar{c}_0^2 c_1 + 2 c_0^3
\bar{c}_0 \bar{c}_1) \psi_0^4 \psi_1 + (3 c_0 \bar{c}_0^2 c_1^2 + 6
c_0^2 \bar{c}_0 c_1 \bar{c}_1 \\
&& + c_0^3 \bar{c}_1^2) \psi_0^3 \psi_1^2 +
(\bar{c}_0^2 c_1^3 + 6 c_0 \bar{c}_0 c_1^2 \bar{c}_1 + 3 c_0^2 c_1 \bar{c}_1^2) \psi_0^2
\psi_1^3 \\
&&+ (2 \bar{c}_0 c_1^3 \bar{c}_1 + 3 c_0 c_1^2 \bar{c}_1^2)
\psi_0 \psi_1^4 + (c_1^3 \bar{c}_1^2) \psi_1^5 +
\mathcal{O} (R).
\end{eqnarray*}

Projecting onto $\psi_0$ and $\psi_1$ respectively and for now
ignoring components with dependence upon $R$ (see \cite{Marzuola}), we
arrive at the finite dimensional
Hamiltonian system of equations given
\begin{eqnarray*}
i \dot{\rho}_0 \langle \psi_0, \psi_0 \rangle & = & - \omega_0 \rho_0
\langle \psi_0, \psi_0 \rangle - \rho_0^3 \bar{\rho}_0^2 \langle \psi_0^5,
\psi_0 \rangle \\
&&- (3 \rho_0 \bar{\rho}_0^2 \rho_1^2 + 6 \rho_0^2 \bar{\rho}_0 \rho_1
\bar{\rho}_1 - \rho_0^3 \bar{\rho}_1^2) \langle \psi_0^4, \psi_1^2 \rangle \\
&& - (2 \bar{\rho}_0 \rho_1^3 \bar{\rho}_1 + 3 \rho_0 \rho_1^2 \bar{\rho}_1^2) \langle
\psi_0^2, \psi_1^4 \rangle ,
\end{eqnarray*}
\begin{eqnarray*}
i \dot{\rho}_1 \langle \psi_1, \psi_1 \rangle & = & - \omega_1 \rho_1
\langle \psi_1, \psi_1 \rangle - \rho_1^3 \bar{\rho}_1^2 \langle \psi_1^5,
\psi_1 \rangle \\
&&- (3 \rho_0^2 \bar{\rho}_0^2 \rho_1 + 2 \rho_0^3 \bar{\rho}_0
\bar{\rho}_1 ) \langle \psi_0^4, \psi_1^2 \rangle \\
&& - (6 \rho_0 \bar{\rho}_0 \rho_1^2 \bar{\rho}_1 + 3 \rho_0^2 \rho_1 \bar{\rho}_1^2
+ \bar{\rho}_0^2 \rho_1^3) \langle
\psi_0^2, \psi_1^4 \rangle
\end{eqnarray*}
and the corresponding conjugate equations.  Note, mass is conserved in
this finite dimensional system, hence we have
\begin{eqnarray*}
|\rho_0|^2 + |\rho_1|^2 = N
\end{eqnarray*}
for all $t$.

Plugging in alternative coordinates designed to give rise to a simple
classification of the finite dimensional dynamics, we set
\begin{eqnarray*}
\rho_0 (t) = A(t) e^{i \theta (t)}
\end{eqnarray*}
and
\begin{eqnarray*}
\rho_1 (t) = (\alpha (t) + i \beta (t) ) e^{ i \theta (t)}.
\end{eqnarray*}
As a result, we have
\begin{eqnarray*}
i \dot{A} - A \dot{\theta} & = & - \omega_0 A - A^5 - 3 A^3 (\alpha +
i \beta)^2 - 6 A^3 (\alpha^2 + \beta^2) \\
&& - 2A (\alpha + i \beta)^2 (\alpha^2 + \beta^2) - A^3 (\alpha - i
\beta)^2 - 3 A (\alpha^2 + \beta^2)^2
\end{eqnarray*}
and
\begin{eqnarray*}
i (\dot{\alpha} + i \dot{\beta}) - (\alpha + i \beta) \dot{\theta} & =
& -
\omega_1 (\alpha + i \beta) - 3 A^2 ( \alpha^2 + \beta^2) (\alpha - i \beta)\\
&& - A^2
(\alpha+i\beta)^3  -2 A^4 (\alpha - i \beta) \\
&&- 6 A^2 ( \alpha^2 + \beta^2) ( \alpha +
i \beta)  - 3 A^4 (\alpha + i \beta) \\
&&+ (\alpha^2+ \beta^2)^2 ( \alpha + i \beta),
\end{eqnarray*}
setting for simplicity $\langle \psi_0^{2j}, \psi_1^{6-2j} \rangle =
1$ for all $j = 0,1,2,3$.  This will simply rescale the dynamical
system and not impact the general shape of the phase diagram.

In the end, we have
\begin{eqnarray*}
\dot{\alpha} & = & (\omega_0  - \omega_1 + 4 A^2 \alpha^2 + 2 ( \alpha^2 +
\beta^2)^2 + 2 (\alpha^2 + \beta^2)( \alpha^2 - \beta^2) ) \beta , \\
\dot{\beta} & = & - (\omega_0 - \omega_1 - 4 A^4 - 4 A^2 \beta^2 + A^2
\alpha^2 + 2 ( \alpha^2 + \beta^2)^2 \\
&&+ 2 (\alpha^2 + \beta^2)( \alpha^2 - \beta^2)) \alpha, \\
\dot{A} & = & - 4A ( A^2 +  (\alpha^2 + \beta^2)) \alpha \beta , \\
\dot{\theta} & = & \omega_0 + A^4 + 10 A^2 \alpha^2  + 2 A^2 \beta^2 +
2 (\alpha^2 + \beta^2)(\alpha^2 - \beta^2) \\
&&+ 3 (\alpha^2 + \beta^2)^2,
\end{eqnarray*}
where
\begin{eqnarray*}
N = A^2 + \alpha^2 + \beta^2.
\end{eqnarray*}

\begin{figure}
\includegraphics[width=1.5in]{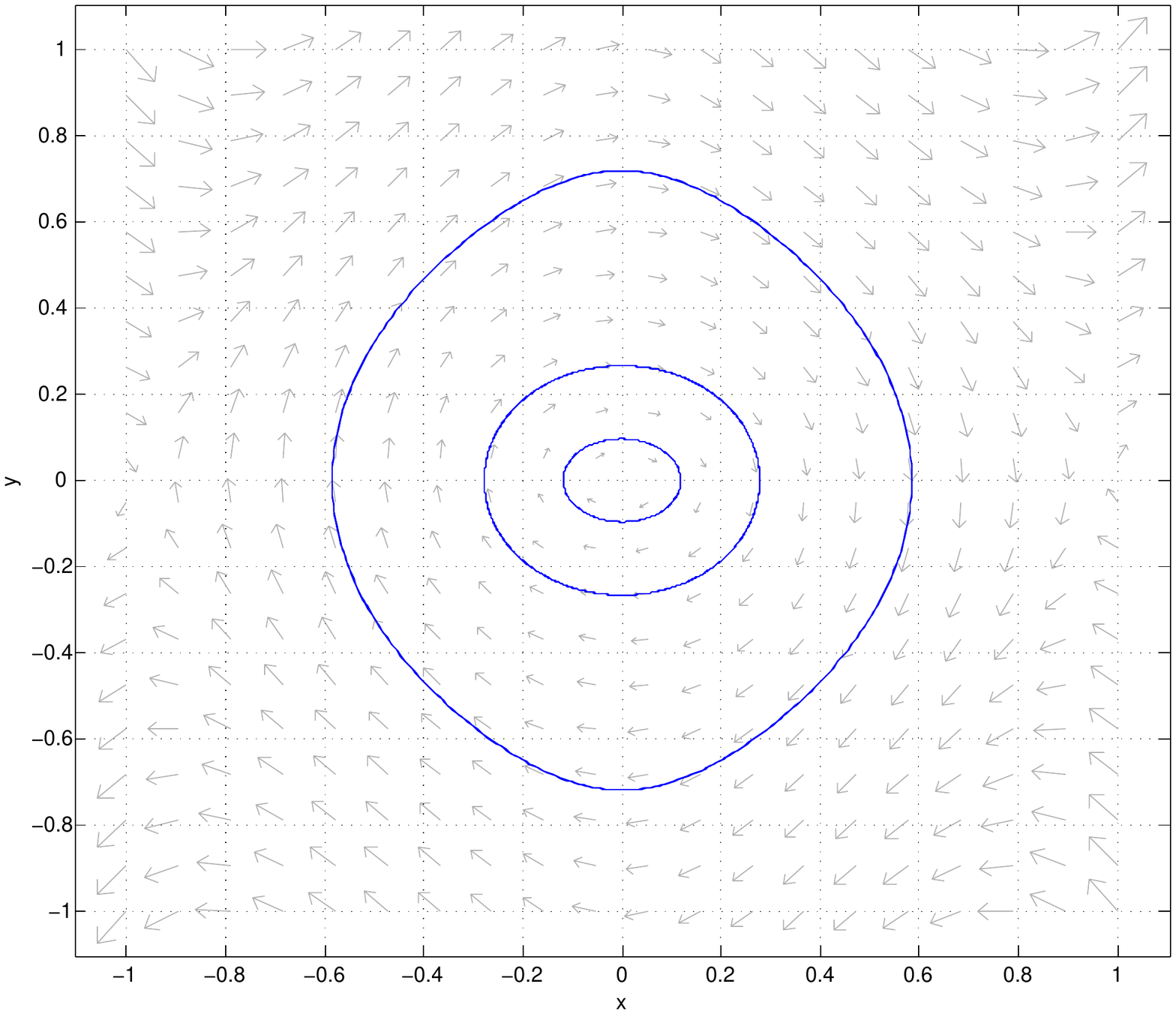}
\includegraphics[width=1.5in]{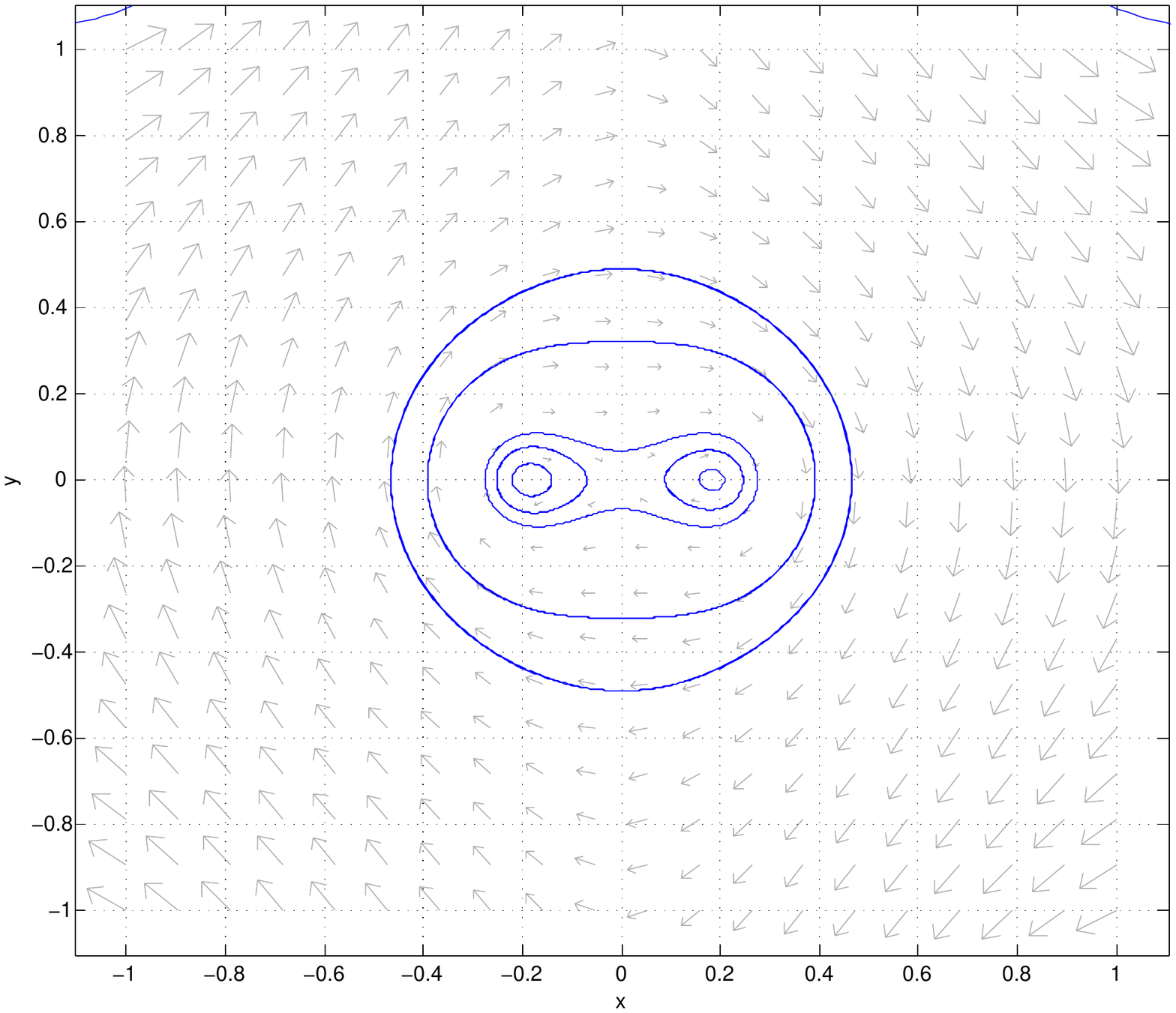}
\includegraphics[width=1.5in]{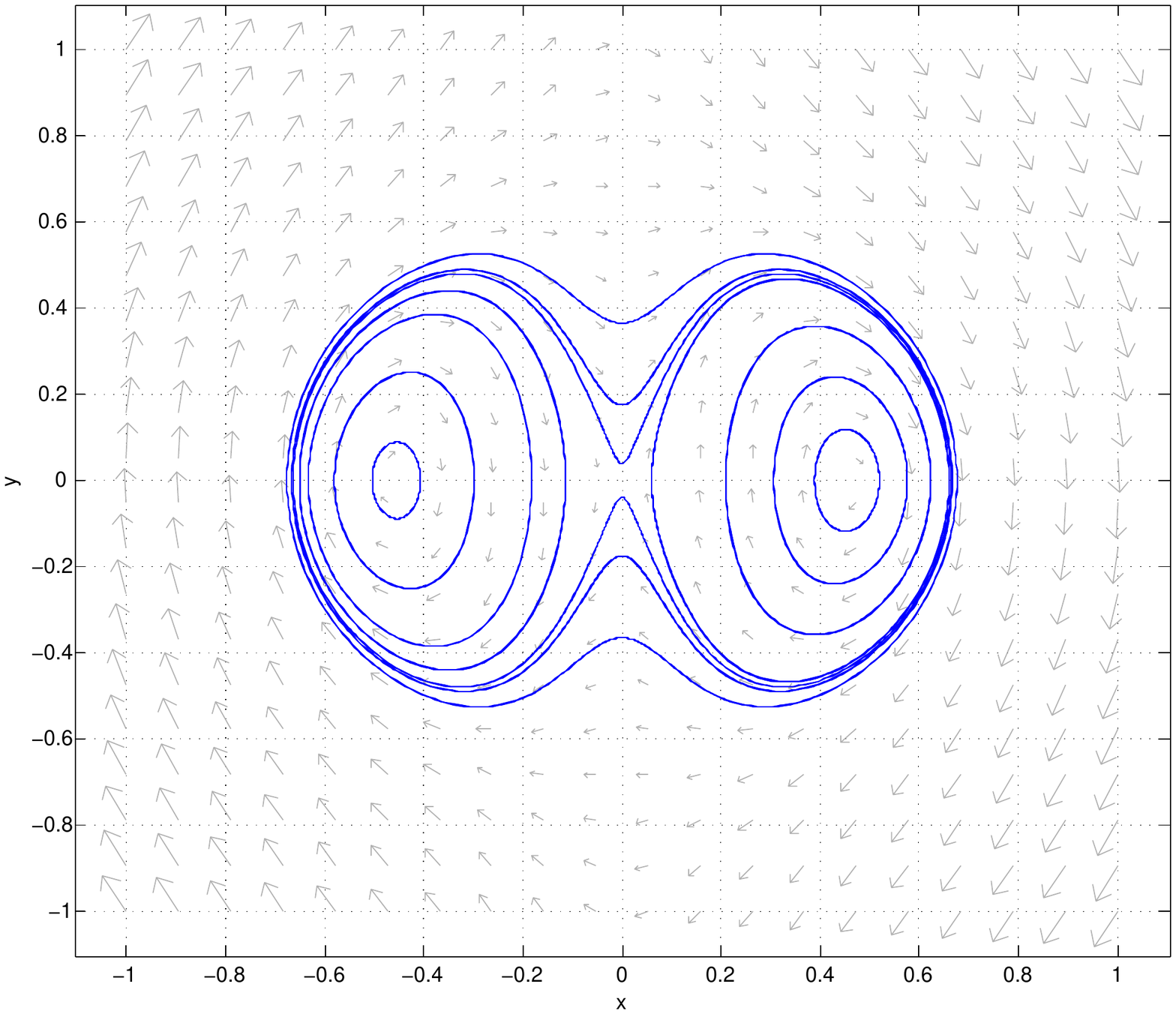}
\caption{Symmetric state for $N-N_{cr}^{FD}<0$ corresponds to an
  equilibrium elliptic point $(\alpha,\beta)=(0,0)$ and asymmetric
  states for $N-N_{cr}^{FD}>0$ correlate to equilibrium points at $(\alpha,\beta)=(\pm \alpha_{cr},0)$.   Plotted phase
  portrait corresponds to  parameter values  $\omega_0 - \omega_1 =
  .1$ and  $N = .1, .2, .5$ respectively.}
\label{fig:pds}
\end{figure}

Using the mass conservation, we can write a closed system for
$(\alpha,\beta)$.  From the equation for $\beta$, it is clear that in
this rescaled dynamical system, we have
\begin{eqnarray*}
N_{cr}^{FD} = \left( \frac{\omega0-\omega1}{4} \right)^{\frac14}.
\end{eqnarray*}
We observe in Figure \ref{fig:pds}, several phase diagrams for varying
values of $N$, which point out the existence of periodic solutions
above, near and below the bifurcation point.  It is our goal in this section purely to give further evidence
that the quintic NLS with double well potential presents similar
dynamics to that of the cubic NLS with double well potential.  For a
dynamics approach to classifying these solutions and studying their
stability properties, we refer to the finite dimensional results
in \cite{Marzuola} for techniques which directly apply to reducible
Hamiltonian systems of this type, particularly for the proof of
existence of periodic orbits and the resulting Floquet stability
analysis.  However, as the intent of this note is to prove asymptotic
stability, we do not explore this topic further here.

\section{Numerical Verification of Hypotheses}
\label{sec:num}

For a potential well
\begin{eqnarray*}
V(x) = \phi_\sigma (x-L) + \phi_\sigma (x+L)
\end{eqnarray*}
with 
\begin{eqnarray*}
\phi_\sigma (x) = \frac{e^{-\frac{x^2}{2*\sigma^2}}}{\sqrt{2 \pi
    \sigma^2}} ,
\end{eqnarray*}
where $\sigma = .001$, $L = 7.5$.  We discretize using a finite element method as in
\cite{Marzuola} on a finer and finer set of grids that are concentrated near the peaks of the delta functions,
we compute using the standard eigenvalue and eigenfunction solvers
from {\it Matlab} the values from \eqref{eq:interaction},
\eqref{eq:interaction2} from hypothesis $(H4)$ are computed as
\begin{eqnarray*}
5 \langle \psi _0^4,  \psi _1 ^2\rangle - \langle \psi _0^6,  1\rangle
\approx .3305
\end{eqnarray*}
and
\begin{eqnarray*}
20  \frac{  5 \langle \psi _0^4,
  \psi _1 ^2\rangle ^2 - \langle \psi _0^6,  1\rangle \langle \psi _0^2,
    \psi _1 ^4 \rangle  }{5 \langle \psi _0^4,  \psi _1 ^2\rangle
    - \langle \psi _0^6,  1\rangle } \!\!\!\!\! &+ & \!\!\!\!\! 160 \frac{\langle \psi _0
^4, \psi _1 ^2\rangle ^2}{\langle \psi _0 ^6,1  \rangle}  -60 \langle
\psi _0 ^2, \psi _1 ^4\rangle \\
&& \approx 9.9143
\end{eqnarray*}
respectively, showing that computationally at least our hypotheses are
valid for a particularly interesting symmetric potential.

\address{ Department of Mathematics,  University of Trieste, Via Valerio  12/1  Trieste, 34127  Italy}
\email{scuccagna@units.it}

\address{Department of Mathematics, University of North Carolina-Chapel Hill \\
Phillips Hall, Chapel Hill, NC  27599, USA}
\email{marzuola@math.unc.edu}


\begin{thebibliography}{CP03}



\bibitem[BP1]{BP1}
V.Buslaev, G.Perelman, {\em Scattering for the nonlinear
Schr\"odinger equation: states close to a soliton\/},  St.
Petersburg Math.J., 4
 (1993), pp.  1111--1142.




\bibitem[BP2]{BP2}
V.Buslaev, G.Perelman, {\em On the stability of solitary waves for
nonlinear Schr\"odinger equations\/},   Nonlinear evolution
equations, editor N.N. Uraltseva, Transl. Ser. 2, 164, Amer. Math.
Soc.,
   pp.  75--98, Amer. Math. Soc., Providence (1995).



\bibitem[CK]{CK}M.Christ, A.Kieslev, Maximal functions
associated with filtrations,  J. Funct. Anal.  179
(200). pp.s 409--425.

 \bibitem[Cu1]{Cu1}
  S.Cuccagna, {\em The Hamiltonian structure of the nonlinear
Schr\"odinger equation and   the  asymptotic stability of its ground
states},      arXiv:0910.3797.

\bibitem[Cu2]{Cu4}   S.Cuccagna, {\em On asymptotic stability in
energy space of  ground states of NLS in 1D},  J.
Differential Equations, 245 (2008), pp.  653-691


\bibitem[Cu3]{Cu5}   S.Cuccagna, {\em A revision of "On asymptotic stability
in energy space of  ground states of NLS in 1D"}, arXiv:0711.4192 .


\bibitem[Cu4]{Cu6}   S.Cuccagna, {\em  Stability of standing waves for NLS with
perturbed Lam\'e potential},  J.
Differential Equations, 223 (2006), pp.  112-160


\bibitem[CM]{cuccagnamizumachi}
S.Cuccagna, T.Mizumachi, {\em On asymptotic stability in energy
space of ground states for Nonlinear Schr\"odinger equations\/},
Comm. Math. Phys., 284
 (2008), pp.  51--87.

 \bibitem[CPV]{CPV} S.Cuccagna, D.Pelinovsky, V.Vougalter, {\em
Spectra of positive and negative energies in the linearization of
the NLS problem},  Comm.  Pure Appl. Math.  {58} (2005), pp. 1--29.


\bibitem[CT]{cuccagnatarulli} S.Cuccagna,  M.Tarulli, {\em On asymptotic stability of standing
waves of discrete
 Schr\"odinger equation  in $  Z$} , SIAM J. Math. Anal.   41,
 (2009), pp. 861-885


\bibitem[DMW]{DMW} V.Duch\^ene, J.L.Marzuola and M.I.Weinstein, {\em Wave
    operator bounds for $1$-dimensional Schr\"odinger operators with
    singular potentials and applications}, to appear in
  J. Math. Phys. (2011).



\bibitem[GS]{GS}
Zhou Gang, I.M.Sigal, {\em Relaxation of Solitons in Nonlinear
  Schr\"odinger Equations with Potential \/}, Advances in
  Math.,
  216 (2007), pp. 443-490.

\bibitem[GSc]{Goldberg}
M.Goldberg, W.Schlag, {\em Dispersive estimates
 for Schr\"odinger   operators in dimensions one and three\/}, Comm. Math. Phys.,
  251 (2004), pp. 157--178.






\bibitem[H]{Har}  E.M. Harrell.  {\em Double Wells}, {  Comm. Math. Phys.}, {  75} (1980), 239-261.

\bibitem[K]{kato}
T.Kato, {\em Wave operators and similarity for some non-selfadjoint
operators \/},  Math. Annalen, 162
 (1966), pp.
  258--269.



\bibitem[KPS]{KPS}   P.G. Kevrekidis, D.E. Pelinovsky, A. Stefanov
{\em Asymptotic stability of small solitons in the discrete
nonlinear Schr\"odinger equation in one dimension} SIAM J. Math.
Anal. 41 (2009),pp. 2010--2030.

\bibitem[KKSW]{kirr} E.Kirr,  P.G.Kevrekidis, E.Shlizerman, M.I.Weinstein.
{\em Symmetry breaking bifurcation in Nonlinear
Schr\"odinger/Gross-Pitaevskii Equations}  SIAM J. Math. Anal. 40
(2008), no. 2, 566--604.

\bibitem[KKP]{KirrKevPel} E.Kirr,  P.G.Kevrekidis, D.E. Pelinovsky.
{\em Symmetry breaking bifurcation in Nonlinear
Schr\"odinger equation with symmetric potentials},  preprint ArXiv:1012.3921 (2010).


\bibitem[KS]{KS} J.Krieger, W.Schlag, {\em Stable manifolds for all
monic supercritical focusing nonlinear Schr\"odinger equations in
one dimension\/}, J. Amer. Math. Soc., 19 (2006), pp. 815--920.



\bibitem[MW]{Marzuola}  J.Marzuola, M.I.Weinstein , {\em Long time dynamics near the symmetry
breaking bifurcation for nonlinear Schr\"odinger / Gross-Pitaevskii
equations }, Discrete and Continuous Dynamical Systems- A, 28
(2010), pp. 1505--1554.





\bibitem[M]{M1}  T.Mizumachi, {\em Asymptotic stability of small
solitons to 1D NLS with potential }, Jour. of Math.   Kyoto
University, 48 (2008), pp. 471-497.



\bibitem[NVT]{NVT} K.Nakanishi, T. Van Phan, T.P. Tsai, {\em Small solutions of nonlinear
Schr\"odinger equations near first excited states}, arXiv:1008.3581 .


\bibitem[PS]{PelinovskyStefanov} D. Pelinovsky, A. Stefanov, {\em
Asymptotic stability of small gap solitons in the nonlinear Dirac
equations},          arXiv:1008.4514.



\bibitem[SW1]{SW1} A.Soffer, M.I.Weinstein, {\em  Multichannel nonlinear
scattering for nonintegrable equations \/}, Comm. Math. Phys., 133
(1990), pp. 116--146






\bibitem[SW2]{SW2}
A.Soffer, M.I.Weinstein, {\em  Multichannel nonlinear scattering II.
The case of anisotropic potentials and data \/},  J. Diff. Eq., 98
 (1992), pp.
  376--390.

\bibitem[SW3]{SW4}    A.Soffer, M.I.Weinstein,
{\em Selection of the ground state for nonlinear Schr\"odinger
equations }, Rev. Math. Phys.  16 (2004), pp.  977--1071.



\bibitem[TY]{TY3}
{  T.P.Tsai, H.T.Yau}, {\em Classification of asymptotic profiles
for nonlinear Schr\"odinger equations with small initial data}, Adv.
Theor. Math. Phys.  {6} (2002), pp.  107--139.





\bibitem[W]{W1}  M.I.Weinstein, {\em Lyapunov stability of ground
states of nonlinear dispersive equations},  Comm. Pure Appl. Math.
 39  (1986), pp.  51--68.













\end{thebibliography}
\end{document}